%%%%%%%%%%%%%%%%%%%%%%%%%%%%%%%%%%%%%%%%%%%%%%%%%%%%%%%%%%%%%%%%%%%%
%%%%                                                            %%%%
%%%%  SEMIPROJECTIVITY FOR CERTAIN PURELY INFINITE C*-ALGEBRAS  %%%%
%%%%  BY JACK SPIELBERG                                         %%%%
%%%%                                                            %%%%
%%%%  PROCESS WITH AMSTEX                                       %%%%
%%%%                                                            %%%%
%%%%%%%%%%%%%%%%%%%%%%%%%%%%%%%%%%%%%%%%%%%%%%%%%%%%%%%%%%%%%%%%%%%%

\input amstex
\documentstyle{amsppt}

%%%%%%%%%%%%%%%%%%%
%%%%%%%%%%%%%%%%%%%
%%%%%%%%%%%%%%%%%%%  MACROS
%%%%%%%%%%%%%%%%%%%
%%%%%%%%%%%%%%%%%%%
\def\qed{\qquad \vrule height6pt width5pt depth0pt}
\def\z{{\bold Z}}

\def\oh{{\Cal O}}

\def\u{{\Cal U}}

\def\t{{\Cal T}}
\def\k{{\Cal K}}
\def\l{{\Cal L}}
\def\real{\mathop{\text{\rm Re}}}

\def\tildee{\widetilde E}

\def\Chi{{\raise2pt\hbox{$\chi$}}}
\def\restrictedto#1{\big|\lower3pt\hbox{$\scriptstyle #1$}}
\def\subsetnoteq{\,\rlap{\raise4true pt\hbox{$\scriptstyle\subset$}}
     \lower3true pt\hbox{$\scriptstyle\not=$}\,}
\def\cstar{$C^*$}

\def\inv{^{-1}}

\def\ofoe{\bigl(o(e)\bigr)}

\def\ofte{\bigl(t(e)\bigr)}

\def\submu{_\mu}
\def\imi{_{i-1}}
\def\ipl{_{i+1}}

\def\jmi{_{j-1}}
\def\jpl{_{j+1}}
\def\jplpl{_{j+2}}
\def\nmi{_{n-1}}
\def\npl{_{n+1}}
\def\nplpl{_{n+2}}
\def\hi{{h_i}}
\def\ki{{k_i}}
\def\tildee{\widetilde{e}}
\def\rank{\hbox{rank}\,}
\def\Xy{\leavevmode
 \hbox{\kern-.1em X\kern-.3em\lower.4ex\hbox{Y\kern-.15em}}}
%%%%%%%%%%%%%%%%%%%
%%%%%%%%%%%%%%%%%%%
%%%%%%%%%%%%%%%%%%%  END MACROS
%%%%%%%%%%%%%%%%%%%
%%%%%%%%%%%%%%%%%%%

\topmatter
\title Semiprojectivity for certain purely infinite \cstar-algebras
\endtitle
\author Jack Spielberg \endauthor
\address Department of Mathematics and Statistics,
Arizona State University,
Tempe, AZ  85287-1804
\endaddress
\email jack.spielberg\@asu.edu\endemail
\abstract
It is proved that classifiable simple separable nuclear purely 
infinite \cstar-algebras having finitely generated K-theory and 
torsion-free K$_1$ are semiprojective.  This is accomplished by 
exhibiting these algebras as \cstar-algebras of infinite directed graphs.
\endabstract
\endtopmatter

\document

\head 1. Introduction \endhead

The notion of semiprojectivity for \cstar-algebras was introduced in 
[EK], with the aim of extending shape theory to the noncommutative 
setting.  An algebra was defined to be semiprojective if it exhibited 
homotopy continuity for maps into an inductive limit \cstar-algebra.  
The definitions of semiprojectivity and shape theory were modified in 
[B1].  For both definitions, the class of algebras  most easily shown 
to be semiprojective were the Cuntz-Krieger algebras ([CK]). In 
effect, the fundamental perturbation lemmas of Glimm ([G]) show that 
the property of being a partial isometry, and the properties of 
equality/orthogonality for two projections, are liftable in the 
appropriate sense.  Since Cuntz-Krieger algebras are defined by 
finitely many generators and relations of these types, it is easy to 
see that they are semiprojective.

The work of Kirchberg, Phillips and R\o rdam has brought attention to 
a class of \cstar-algebras defined by the essential features of the 
Cuntz-Krieger algebras; namely, the simple separable nuclear purely 
infinite \cstar-algebras.  In the presence of the universal 
coefficient theorem, it was shown that these algebras 
were classified (up 
to strong Morita equivalence) by K-theory ([K], [P]), and that the 
K-theory can be any pair of countable abelian groups.  It is natural 
to ask whether these algebras also are semiprojective.  
Blackadar has shown
that finitely generated K-theory is necessary for 
semiprojectivity ([B2, Corollary 2.11]).  In the same paper he 
established 
semiprojectivity  in the case where the $K_0$-group is free and the 
$K_1$-group is trivial.  His arguments were generalized by Szymanski 
([Sz1]) to the case where the $K_0$-group is finitely generated and 
the $K_1$-group is free with $\rank K_1\le\rank K_0$.
In this paper, we establish semiprojectivity without the restriction 
on the ranks of the $K$-groups.  In particular we show that the 
algebra $P_\infty$ having $K_0=0$ and $K_1=\z$ is semiprojective.  
Although our proof also relies on realizing the algebras as 
\cstar-algebras of infinite graphs, our techniques are substantially 
different from those in [B2] and [Sz1].  We wish to draw attention to 
our Lemma 3.6, which we feel is a novel addition to the field of 
noncommutative general topology.

For (irreducible) finite graphs, K$_0$ and K$_1$ have equal rank.  To 
model the general situation we must use infinite graphs.  In another 
paper ([Sp]) we proved that a graph \cstar-algebra has an inductive 
limit decomposition over the (directed set of the) finite subgraphs.
(Similar results have been obtained by different methods in [RS].)
This provides a strategy for constructing infinite graphs whose 
\cstar-algebra has specified K-theory, which we follow in section 2 
of this paper.  In particular, we provide an elegant graph whose 
\cstar-algebra is the classifiable example having $K_0=(0)$ and 
$K_1\cong\z$.  (A row-finite graph with the same \cstar-algebra was 
given by Neub\"{u}ser in [N], but has not been found useful for proving 
semiprojectivity.  
It is also shown in [N] that the $K$-theory of irreducible Exel-Laca
\cstar-algebras (of infinite matrices; see [EL]) can be any pair of 
countable abelian groups with free $K_1$.  In [Sz2] Szymanski proved 
the analogous result for \cstar-algebras of row-finite graphs.)
In section 3 we again use the inductive limit 
decomposition to prove semiprojectivity.
\smallskip
The figures in this paper were prepared with \Xy-pic.

\medskip
We recall here the definition of semiprojectivity from [B1] (slightly 
modified).

\definition{Definition} The \cstar-algebra $A$ is called {\it 
semiprojective} if the following holds. If $B$ is
 a \cstar-algebra containing
 a directed  
family $\l$ of ideals with closure $I$, and  $B/I$ is 
isomorphic to $A$, then there is an ideal $J\in\l$ and a 
$*$-homomorphism from $A$ to $B/J$ lifting the quotient map.
\enddefinition

Let $E$ be a directed graph with vertices $E^0$ and (oriented)
edges $E^1$.  We will follow the notations in [Sp], using $o$ and $t$ 
for the origin and terminus of an edge. 
However, we will omit the subscript $+$ indicating
(positively) oriented edges, 
as in this paper we will not make use of the negatively oriented ones.
We will let the edges and vertices also denote the generating elements 
of the algebra $\oh(E)$. Potential confusion will
be avoidable from the context.
We will let $D$ denote the set of vertices of $E$ having infinite exit 
valence:
$$D:=\{u\in E^0\bigm|o\inv(u)\hbox{\ is\ infinite}\}.$$
We recall from [Sp, Theorem 2.21], the relations defining $\oh(E)$.  
We  refer to these 
as condition $(\oh)$.

\definition{Condition $(\oh)$}
\smallskip\noindent
\item{(o1)} $\{e\mid e\in E^1\}$ are partial isometries, and $\{u\mid 
u\in E^0\}$ are projections.
\item{(o2)} $uv=0$ if $u\not= v$, for all $u$, $v\in E^0$.
\item{(o3)} $e^*e=t(e)$, for all $e\in E^1$.
\item{(o4)} $o(e)\cdot e=e$, if $o(e)\in D$.
\item{(o5)} If $e\not=f$ and $o(e)=o(f)\in D$, then $e^*f=0$.
\item{(o6)} If $u\not\in D$ then $u=\sum\{ee^*\mid e\in o\inv(u)\}$.
\enddefinition

We will be primarily interested in graphs with special properties.  We 
set out these properties here.  We consider an infinite directed 
graph $E$ having a chain $F_0\subseteq F_1\subseteq 
\ldots\subseteq E$ of finite subgraphs with $\cup F_n = E$. Let $D$ 
be as above.  We recall from [Sp, Theorem 2.34], that the inclusion of graphs 
$F_n\subseteq E$ determines a Toeplitz graph algebra $\t\oh(F_n,S_n)$, 
where $S_n=\{u\in F_n^0\mid o\inv(u)\cap E^1\subseteq F_n^1\}$.  We 
will let $\t\oh(F_n)$ denote this algebra.
It follows from [Sp, Theorem 2.20], that the  
subalgebra $\t\oh(F_n)$ is generated by the edges and vertices of 
$F_n$ subject to the same relations $\oh$, with the exception that 
the set $D$ above must be enlarged to include $F_n^0\setminus S_n$.

\remark{Remarks 1.1}  We recall from [Sp, Theorem 2.35], that 
$$\oh(E)=\overline{\cup_n\t\oh(F_n)}.$$
We recall from [Sp, Theorem 3.15], that if $E$ is an irreducible graph 
that is not a cycle, then $\oh(E)$ is simple, nuclear, and purely 
infinite.
\endremark

\definition{Condition (a)}
\item{(a1)} $D\subseteq F_0^0$
\item{(a2)} $F_n$ is irreducible and not a cycle, $n\ge0$
\item{(a3)} $o\inv(F_n^0)\cap t\inv(F_n^0)\cap E^1 \subseteq 
F_n^1$, $n\ge0$
\item{(a4)} $o\inv(F_n^0\setminus D)\subseteq F_n^1$, 
$n\ge0$
\item{(a5)} $t\inv(F_n^0)\subseteq F_{n+1}^0$, $n\ge0$.
\enddefinition

We may loosely summarize condition (a) in the following way:  
$E\setminus D$ is locally finite and flows towards $D$, while $D$ 
points to all of the subgraphs.

\definition{Condition (b)}
\item{(b1)} $K_*\bigl(\oh(E)\bigr)$ is finitely generated
\item{(b2)} For each $u\in E^0\setminus D$ there exist $v\in 
F^0_0\setminus D$ 
and a graph isomorphism $\gamma:o\inv(u)\to o\inv(v)$ such that 
$\bigl[t(e)\bigr]=\bigl[t(\gamma(e))\bigr]$ in $K_0\bigl(\oh(E)\bigr)$ 
for $e\in o\inv(u)$.
\enddefinition

Condition (b2) is a technical aid in the proof of Theorem 3.1, and it 
is evident in the examples to which we will apply that theorem.

\remark{Remarks 1.2}
(i) Since all $F_n$ are irreducible, $E$ is 
irreducible.
\smallskip\item{}(ii) Since $F_0$ is finite, (a1) implies that $D$ is finite.
\smallskip\item{}(iii) Assuming $F_{n+1}\not=F_n$ for all $n\ge0$, it follows 
from (a3) that $F_n^0\setminus F_{n-1}^0\not=\emptyset$.  Since $F_n$ 
is irreducible, then
$$t\inv(F_n^0\setminus F_{n-1}^0)\cap o\inv(F_{n-1}^0) 
\not=\emptyset.$$
\par By (a4) it follows that
$$t\inv(F_n^0\setminus F_{n-1}^0)\cap o\inv(D) \not=\emptyset.$$
\par In particular, $D\not=\emptyset$.
\smallskip\item{}(iv) If $\{F_n\mid n=0,1,2,\ldots\}$ is replaced by a 
subsequence, properties (a1) -- (a5) still hold.
\smallskip\item{}(v) By passing to a subsequence if necessary, we may assume 
that
\smallskip
\item{(a6)} $t\inv(F_n^0\setminus F_{n-1}^0)\cap 
o\inv(u)\not=\emptyset$, $u\in D$, $n\ge0$.
\smallskip
\smallskip
\item{}(vi) It follows from (b1)
and Remark 1.1 that by enlarging $F_0$ 
we may assume
\smallskip
\item{(b3)} $i_*:K_1\bigl(\t\oh(F_0)\bigr)\to 
K_1\bigl(\oh(E)\bigr)$ is onto.
\smallskip
\item{}Since $\oh(E)$ is simple and purely infinite, it follows from
[C2] that by enlarging $F_1$ we may assume
\smallskip
\item{(b4)} If $v_1,v_2\in F_0^0$ satisfy $[v_1]=[v_2]$ in 
$K_0\bigl(\oh(E)\bigr)$, then $v_1\sim v_2$ in $\t\oh(F_1)$.
\endremark

\remark{Remark 1.3} Because of (a1), the relations defining 
representations of $\t\oh(F_n)$ are exactly the same as $(\oh)$ above 
(when applied to $F_n$).  (In this case, the representation is faithful 
if and only if for each $u\in D$, $u\not=\sum\{ee^*\mid e\in 
o\inv(u)\}$.  See [Sp, Theorem 2.20].)
\endremark

\head 2. K-theory \endhead

In Theorem 2.3 below we construct graphs whose algebras have certain 
specified K-theory.  To this end we must first compute the K-theory of 
the Toeplitz graph algebras $\t\oh(F_n)$, and also study the maps in 
K-theory induced by the inclusions $\t\oh(F_n)\subseteq\t\oh(F\npl)$.

Let $F$ be a finite  graph satisfying
 condition (I) of [C1].  Let $S_F\subseteq F^0$, and let 
$\t\oh(F)$ denote $\t\oh(F,S_F)$.  We recall from [Sp, Theorem 2.33],
that there is a 
short exact sequence
$$0\longrightarrow I_F\longrightarrow\t\oh(F)
\buildrel\pi\over\longrightarrow \oh(F) 
\longrightarrow0,$$
where $I_F\cong\bigoplus_{S_F^c}\k$.  We obtain a long exact sequence 
in K-theory:
$$0\longrightarrow K_1\t\oh(F)\longrightarrow K_1\oh(F)
\buildrel \partial\over\longrightarrow \bigoplus_{S_F^c}\z 
\longrightarrow K_0\t\oh(F) \longrightarrow K_0\oh(F)\longrightarrow 
0.$$
From [C1] we have $K_*\oh(F)$:
$$\eqalign{K_0\oh(F)&= \z^{F^0}\Big/\Bigl\{[u]=\sum_{e\in 
o\inv(u)}[t(e)]\Bigr\}\cr
\noalign{\smallskip}
K_1\oh(F)&= \bigl\{x\in\z^{F^0}\bigm| 
x(u)=\sum_{e\in t\inv(u)}x\ofoe\bigr\}.\cr}$$
For $x\in K_1\oh(F)$ we will further define $x$ on $F^1$ by setting
$$x(e):=x\ofoe.$$

\proclaim{Lemma 2.1}  If $x\in K_1\oh(F)$, then $(\partial x)_u = 
-\sum\bigl\{x(e)\bigm|e\in t\inv(u)\bigr\}$,
for $u\in S_F^c$.
\endproclaim

\demo{Proof} We will let $e\in F^1$ represent the corresponding generator of 
$\t\oh(F)$, and $\tildee$ represent the generator 
$\pi(e)\in\oh(F)$.
For $u\in F^0$ let $\varepsilon_u\in\t\oh(F)$ be given by
$$\varepsilon_u= u\;-\sum_{e\in o\inv(u)}ee^*.$$
Then $\varepsilon_u\not=0$ if and only if $u\not\in S_F$, and 
$\pi_*(\varepsilon_u)=0$ in $\oh(F)$ for all $u$.  Moreover, the classes 
$\bigl\{[\varepsilon_u]\bigm|u\in S_F^c\bigr\}$ 
form a basis for $K_0(I_F)$.

Now let $x\in K_1\oh(F)$.  We construct a unitary representing $x$ as 
in [R, section 2].  With $|x|=\sum_{e\in F^1}|x(e)|$, 
we let $V(x)\in 
M_{|x|}\bigl(\oh(F)\bigr)$ be a diagonal matrix with $|x(e)|$ entries 
equal to $\tildee$, if $x(e)>0$, and with $|x(e)|$ entries 
equal to $\tildee^*$, if $x(e)<0$.  Then $V(x)^*V(x)$ and 
$V(x)V(x)^*\in M_{|x|}\bigl(\hbox{span}\,\{\tildee\tildee^*\bigm| e\in 
F^1\}\bigr)$, and $1-V(x)^*V(x)$ is equivalent to $1-V(x)V(x)^*$.  
Let $W(x)\in M_{|x|}\bigl(\hbox{span}\,\{\tildee\tildee^*\bigm| e\in 
F^1\}\bigr)$ with $W(x)^*W(x)=1-V(x)^*V(x)$ and $W(x)W(x)^* = 
1-V(x)V(x)^*$.  Then $U(x)=V(x)+W(x)\in M_{|x|}\bigl(\oh(F)\bigr)$ is 
unitary, and $[U(x)]=x$.

Next note that if we replace $\tildee$ by $e$ throughout, we obtain 
partial isometries $U_0(x)$, $V_0(x)$, $W_0(x)\in M_{|x|}\bigl( 
\t\oh(F)\bigr)$ that are lifts (modulo $M_{|x|}(I_F)$) of $U(x)$, 
$V(x)$, $W(x)$.  Therefore,
$$\partial x=\bigl[I-U_0(x)^*U_0(x)\bigr] - 
\bigl[I-U_0(x)U_0(x)^*\bigr].$$
The fact that $\pi_*\circ\partial(x)=0$ accounts for all terms in the 
above involving the $ee^*$.  What remains are the occurrences of 
the $\varepsilon_u$.  We obtain
$$\eqalign{\partial x&=\left(\sum_{x(e)>0}|x(e)|\sum_{u\not=t(e)} 
[\varepsilon_u]\quad + \quad\sum_{x(e)<0}|x(e)|\sum_u[\varepsilon_u]\right)\cr
&\qquad\qquad - \left(
\sum_{x(e)>0}|x(e)|\sum_u
[\varepsilon_u]\quad +\quad \sum_{x(e)<0}|x(e)|\sum_{u\not=t(e)}
[\varepsilon_u]\right)\cr
&=\sum_ex(e)\left(\sum_{u\not=t(e)}[\varepsilon_u]\quad-\quad\sum_u 
[\varepsilon_u]\right)\cr
&=-\sum_e x(e)[\varepsilon_{t(e)}]\cr
&=-\sum_{u\in S_F^c}\left(\sum_{e\in t\inv(u)}x(e)\right) 
[\varepsilon_u]. \qed\cr}$$
\enddemo

Thus we find that
$$\eqalignno{
K_1\bigl(\t\oh(F)\bigr)\cong\ker(\partial)&= \left\{x\in\z^{F^0} 
\biggm| \sum_{e\in t\inv(u)}x(e)=0,\ u\in S_F^c\right\},\cr
\hbox{Im}\,(\partial)&=\left\{\biggl(\sum_{e\in t\inv(u)}x(e)
\biggr)_{u\in S_F^c}\biggm|x\in 
K_1\bigl(\oh(F)\bigr)\right\}.\cr}$$

Now let $F\subseteq G$ be an inclusion of finite directed graphs 
satisfying condition (I).  We assume that $S_F\subseteq S_G$ (as would 
be the case if $G\subseteq E$, and $S_F$, $S_G$ are defined relative 
to $E$).  Let $i:\t\oh(F)\subseteq\t\oh(G)$ denote the inclusion.

\proclaim{Lemma 2.2}  Let $x\in K_1\bigl(\t\oh(F)\bigr)$.  Then
$$(i_*x)_u=\cases x(u),&\text{if $u\in F^0$,}\cr
0,&\text{if $u\in G^0\setminus F^0$.}\cr \endcases$$
\endproclaim

\demo{Proof} When viewed as an element of $K_1\bigl(\oh(F)\bigr)$, $x$ is 
represented by the diagonal unitary matrix $U(x)\in 
M_{|x|}\bigl(\oh(F)\bigr)$, as in the proof of Lemma 2.1.  Then 
$i_*(x)$ is represented by $U(x)+I-U(x)^*U(x)\in M_{|x|}\bigl( 
\oh(G)\bigr)$.  But this is exactly the unitary representing the 
element of $\z^{G^0}$ given in the statement of the Lemma. \qed
\enddemo

Let $j_F:I_F\subseteq \t\oh(F)$ be the inclusion.
Then $i_* \circ j_{F\,*}:K_0(I_F)\to K_0\bigl(\t\oh(G)\bigr)$ is 
easily computed to be
$$i_* \circ j_{F\,*}[\varepsilon_{u,F}]= [\varepsilon_{u,G}] + \sum\Bigl\{
[ee^*]\bigm|
e\in o\inv(u)\cap G^1\setminus F^1\Bigr\}.$$
We omit the proof.

\proclaim{Theorem 2.3} Let $G_0$ and $G_1$ be finitely 
generated abelian groups 
with $G_1$ free.  Then there is an irreducible non-cycle graph $E$ 
with $K_*\oh(E)\cong(G_0,G_1)$.  If $\rank G_0\not=\rank G_1$, then 
$E$ can be chosen along with a chain of subgraphs $(F_n)_{n\ge0}$ so 
that conditions (a) and (b) are satisfied.
\endproclaim

\demo{Proof} {\it Case (i):\/}  $\rank G_0=\rank G_1$.  The existence of such 
a graph $E$ is known, e.g. from [R].  Since we will need a specific graph 
with this property below, we describe it here.  Suppose that 
$$\eqalign{G_0&=\z^\ell\oplus \z/(n_1)\oplus \ldots \oplus \z/(n_k)\cr
G_1&=\z^\ell,\cr}$$
where $k$, $\ell\ge0$, and $n_i\ge2$.  We construct a graph $F$ with 
$k+\ell+1$ vertices $v_1$, $\ldots$, $v_\ell$, $z$, $w_1$, $\ldots$, 
$w_k$, by including edges as follows:
$$\eqalign{v_i&\longrightarrow\{v_1,\ldots,v_\ell,z\}\cr
z&\buildrel n_1\over\longrightarrow w_1\cr
w_i&\buildrel n_{i+1}\over\longrightarrow w_{i+1},\quad 1\le i<k,\cr
w_k&\longrightarrow\{v_1,\ldots,v_\ell,z\},\cr}$$
and an additional loop at each vertex.  It is clear that the graph 
$F$ is irreducible (and even primitive).  Thus, with the vertices 
used in the order listed above, the incidence matrix of the graph is
$$ \setbox0=\vbox{\def\cr{\crcr\noalign{\kern2pt\global\let\cr=\endline}}
 \offinterlineskip
  \halign{$#$\hfil\kern2pt\kern\wd0&\thinspace\hfil$#$\hfil
    &\quad\vrule#&\quad\hfil$#$\hfil\crcr
        \strut&\;\ell+1&\omit&k\cr 
        \strut \ell\;&
	  \matrix \cr &\hbox{1's}\cr \cr\endmatrix
	&&
	  \matrix \cr &0&\cr\cr\endmatrix
	\cr&\multispan3\hrulefill\cr
        \vphantom{o}&&\cr
        \strut k\;&
	  \matrix &0\cr\endmatrix
	&&
	  \matrix n_1\cr &\ddots\cr &&n_k\cr\endmatrix
	\cr\vphantom{o}&&\cr
        &\multispan3\hrulefill\cr
        \vphantom{o}&&\cr
        \strut 1\;&
	  \matrix &\hbox{1's}\cr\endmatrix
	&&
	  \matrix &0&\cr\endmatrix
	\cr\crcr\omit\strut\cr}}
  \setbox2=\vbox{\unvcopy0 \global\setbox1=\lastbox}
  \setbox2=\hbox{\unhbox1 \unskip \global\setbox1=\lastbox}
  \setbox2=\hbox{$\kern\wd1\kern-\wd0 \left( \kern-\wd1
   \global\setbox1=\vbox{\box1\kern2pt}
   \vcenter{\kern-\ht1 \unvbox0 \kern-\baselineskip} \,\right)$}
  \null\;\vbox{\kern\ht1\box2}\quad+\quad I.$$
(If $k=0$ we use the sum of an $(\ell+1)\times(\ell+1)$ matrix of 1's with 
the identity.)  If $x\in \z^{F^0}$ represents an element of 
$K_1\oh(F)$, we obtain the following relations:
$$\leqalignno{\sum_{i=1}^\ell x(v_i)+x(w_k)&=0,\hbox{\ from\ }v_i 
\hbox{\ or\ }z,&\hbox{(k1)}\cr
x(z)&=0,\hbox{\ from\ }w_1,&\hbox{(k2)}\cr
x(w_i)&=0,\;i<k,\hbox{\ from\ }w\ipl.\cr}$$
It follows that $K_1\oh(F)\cong\z^\ell$, with basis $x_1$, $\ldots$, 
$x_\ell$, where
$$x_i(s)=\cases \phantom{-}1,&\text{if $s=v_i$}\cr -1,&\text{if  
$s=w_k$}\cr\phantom{-} 0,&\text{else.}\endcases$$

$K_0\oh(F)$ is generated by the classes of the vertices, with relations
$$\leqalignno{\sum_{i=1}^\ell[v_i] + [z]&=0,\hbox{\ from\ }v_i 
\hbox{\ or\ }w_k,&\hbox{(k3)}\cr
n_i[w_i]&=0,\;1\le i\le k,\hbox{\ from\ }z,w_1,\ldots,w_{k-1}.&
\hbox{(k4)}\cr}$$
It follows that $K_0\oh(F)\cong\z^\ell\oplus
\z/(n_1)\oplus \ldots \oplus \z/(n_k)$, with the cyclic subgroups 
generated by $[v_1]$, $\ldots$, $[v_\ell]$, $[w_1]$, $\ldots$, $[w_k]$.

\smallskip\noindent
{\it Case (ii):\/}  $\rank G_0<\rank G_1$.  Let $G_1\cong\z^{\ell-1}$ 
and $G_0\cong\z^{\ell-p-1}\oplus
\z/(n_1)\oplus \ldots \oplus \z/(n_k)$, where $0<p<\ell$ and $\ell>1$. 
We build a larger graph, $E$, from $F$ as follows.  Include a vertex 
$u$, and for $1\le i\le p+1$ include the graph shown in Figure 1.

\bigskip
\input epsf         % defines \epsfbox and supporting macros
\epsfverbosetrue    % messages will show height and width
\centerline{\epsfbox{fig1.epsi}}
\bigskip

For $n\ge0$ we let $F_n$ be the subgraph of $E$ with
$$\eqalign{F_n^0&=F^0\cup\{u\}\cup\bigl\{a_{ij},\,b_{ij}\bigm|1\le i\le 
p+1,\  
1\le j\le n+1\bigr\}\cr
F_n^1&=o\inv(F_n^0)\cap t\inv(F_n^0).\cr}$$
We will let $\t\oh(F_n)$ denote 
$\t\oh\bigl(F_n,F_n^0\setminus\{u\}\bigr)$.
We compute $K_*\oh(F_n)$.  Let $x\in\z^{F_n^0}$ represent an element 
of $K_1\oh(F_n)$.  Note that the relations (k1) and (k2) still hold 
((k1) from the vertex $z$).  For $1\le i\le p+1$, consideration of the 
vertex $v_i$ gives
$$\sum_{j=1}^\ell x(v_j) + x(w_k) + x(a_{i1})=0.\leqno\hbox{(k5)}$$
Combined with (k1), we find that
$$x(a_{i1})=0,\ 1\le i\le p+1.\leqno\hbox{(k6)}$$
Consideration of the vertex $u$ gives
$$x(u)=\sum_{i=1}^{p+1}x(v_i).\leqno\hbox{(k7)}$$
Finally, consideration of the remaining vertices show that
$$\eqalign{x(a_{ij})&=0,\cr x(b_{ij})&=-x(u),\cr}\quad 1\le i\le p+1,\ 1\le 
j\le n+1.\leqno\hbox{(k8)}$$
It follows from (k1), (k2), (k6), (k7), and (k8) that 
$K_1\oh(F_n)\cong\z^\ell$, with the same generators as for 
$K_1\oh(F)$.  Hence (k7) implies that 
$K_1\t\oh(F_n)\cong\z^{\ell-1}$.  Moreover it is clear that the 
inclusion $\t\oh(F_n)\subseteq \t\oh(F\npl)$ induces an isomorphism 
in $K_1$.  Hence $K_1\oh(E)\cong\z^{\ell-1}$.

It also follows from (k7) that $\partial$ is onto, so that 
$K_0\t\oh(F_n)\cong K_0\oh(F_n)$.  To compute $K_0\oh(F_n)$, notice 
that (k3) and (k4) still hold ((k3) from the vertex $w_k$; or if 
$k=0$, from the vertex $z$).  For $1\le j\le p+1$, consideration of 
the vertex $v_i$ gives
$$[u]+\sum_{j=1}^\ell[v_j] + [z]=0,$$
and hence
$$[u]=0.\leqno\hbox{(k9)}$$
Consideration of the vertices $a_{ij}$ shows that
$$\eqalign{[v_i]&=0,\quad 1\le i\le p+1\cr
[a_{ij}]&=0,\quad 1\le i\le p+1,\ 0\le j\le n.\cr} \leqno\hbox{(k10)}$$
Finally, consideration of the vertices $b_{ij}$ shows that
$[b_{ij}] = 2[b_{ij}]+[a_{ij}]$, and hence
$$[b_{ij}]=-[a_{ij}],\quad 1\le i\le p+1,\ 0\le j\le n+1. 
\leqno\hbox{(k11)}$$
It follows from (k3), (k4), (k9), (k10), and (k11) that 
$K_0\oh(F_n)\cong
\z^\ell\oplus \z/(n_1)\oplus \ldots \oplus \z/(n_k)$, with basis for 
the free summand given by
$$[v_{p+2}], \ldots,[v_\ell],\;[a_{1,n+1}],\ldots,[a_{p+1,n+1}].$$
Thus the map $i_*:K_0\t\oh(F_n)\to K_0\t\oh(F_{n+1}$ kills the last 
$p+1$ of the basis elements.  Thus $K_0\oh(E) = 
\lim\limits_\to K_0\t\oh(F_n)\cong
\z^{\ell-p-1}\oplus \z/(n_1)\oplus \ldots \oplus \z/(n_k)$.

\smallskip\noindent
{\it Case (iii):\/}  $\rank G_0>\rank G_1$.  Let $G_1\cong\z^{\ell-p}$ 
and $G_0\cong\z^{\ell}\oplus
\z/(n_1)\oplus \ldots \oplus \z/(n_k)$, where $0<p\le\ell$ and 
$\ell\ge1$. 
We build a larger graph, $E$, from $F$  by attaching to each $v_i$,
$1\le i\le p$,  the graph shown in Figure 2.

\bigskip
\input epsf         % defines \epsfbox and supporting macros
\epsfverbosetrue    % messages will show height and width
\centerline{\epsfbox{fig2.epsi}}
\bigskip

For $n\ge0$ we let $F_n$ be the subgraph of $E$ given by
$$\eqalign{F_n^0&=F^0\cup\{u_i\mid1\le i\le p\}\cup
\bigl\{a_{ij},\,b_{ij}\bigm|1\le i\le p,\ 
1\le j\le n+2\bigr\}\cr
F_n^1&=o\inv(F_n^0)\cap t\inv(F_n^0),\cr}$$
and we let let $\t\oh(F_n)$ denote 
$\t\oh\bigl(F_n,F_n^0\setminus\{u_i\mid 1\le i\le p\}\bigr)$.
We first compute $K_1\oh(F_n)$.  Let $x\in\z^{F_n^0}$ represent an 
element of $K_1\oh(F_n)$.  We still have (k1) and (k2), as in Case 
(ii).  Consideration of vertex $a_{ij}$ gives
$$x(a_{ij})+x(b_{ij})=0,\quad 1\le i\le p,\ 1\le j\le n+2.\leqno 
\hbox{(k12)}$$
Consideration of vertex $v_i$, $1\le i\le p$, gives
$$\sum_{j=1}^\ell x(v_j) + x(w_k) + x(a_{ij}) + x(b_{ij})=0,$$
which is already contained in (k1) and (k12).  Consideration of 
vertex $u_i$ gives
$$x(u_i)=x(v_i),\quad1\le i\le p.\leqno\hbox{(k13)}$$
Finally, consideration of vertex $b_{ij}$ gives
$$x(b_{ij})=x(u_i),\quad1\le i\le p,\ 0\le j\le n+2, 
\leqno\hbox{(k14)}$$
(with a slightly different computation in the case $j=n+2$).
It follows from (k3), (k4), (k12), (k13), and (k14)
 that $K_1\oh(F_n)\cong\z^\ell$. Moreover (k13) implies that 
 $\partial$ is onto, so that $K_1\t\oh(F_n)\cong\z^{\ell-p}$, with 
 basis $x_{p+1}$, $\ldots$, $x_\ell$ (where the $x_i$ are as defined 
 in Case (ii)).  It follows that $\t\oh(F_n)\subseteq\t\oh(F_{n+1}$ 
 induces an isomorphism in $K_1$, so that $K_1\oh(E)\cong\z^{\ell-p}$.

Since $\partial$ is onto, 
$K_0\t\oh(F_n)\cong K_0\oh(F_n)$.  To compute $K_0\oh(F_n)$, notice 
that (k3) and (k4) still hold, as in Case (ii).  Consideration of 
vertex $v_i$ gives
$$\sum_{j=1}^\ell[v_j] + [z] + [u_i]=0,$$
and hence, from (k3),
$$[u_i]=0,\quad1\le i\le p.\leqno\hbox{(k15)}$$
Consideration of vertex $a_{ij}$ gives
$$\eqalign{[a_{i1}]&=-[v_i],\quad 1\le i\le p\cr
[a_{ij}]&=-[b_{i,j-1}],\quad 1\le i\le p,\ 1\le j\le n+2.\cr} 
\leqno\hbox{(k16)}$$
Consideration of the vertices $b_{ij}$ inductively in $j$ gives
$$[b_{ij}]=0,\quad1\le i\le p,\ 1\le j\le n+2.\leqno\hbox{(k17)}$$
It follows from (k3), (k4), (k15), (k16), and (k17) that 
$K_0\oh(F_n)\cong K_0\oh(F)$, with the same generators and relations.  
Hence $\t\oh(F_n)\subseteq\t\oh(F_{n+1})$ induces an isomorphism in 
$K_0$, so that $K_0\oh(E)\cong K_0\oh(F)$.

Finally, we claim that in Cases (ii) and (iii) the graph $E$ and chain 
of subgraphs $(F_n)_{n\ge0}$ satisfy conditions (a) and (b).  
Conditions (a) and (b1) are clear.  For (b2), in Case (ii) we 
identify $a_{ij}$ and $b_{ij}$ with $a_{i1}$ and $b_{i1}$, 
respectively, in the obvious way, for $j\ge 2$.  In Case (iii) we  
identify $a_{ij}$ and $b_{ij}$ with $a_{i2}$ and $b_{i2}$, 
respectively, in the obvious way, for $j\ge 3$. \qed
\enddemo

\remark{Remark 2.4}  In the case where $G_0$ is also torsion-free, the 
graph $E$ can be simplified.  We show in Figure 3 a graph $E$ for 
which $K_0\oh(E)\cong0$ and $K_1\oh(E)\cong\z$.  (Neub\"{u}ser,
in [N, Proposition 3.17],
has given 
a row-finite graph having the same \cstar-algebra.)
If the graph is 
expanded to have $\ell+1$ `columns' joined at the one vertex having 
infinite valence, a graph is obtained for which $K_*\cong(0,\z^\ell)$.  We 
leave it to the interested reader to check the details.
\endremark

\bigskip
\input epsf         % defines \epsfbox and supporting macros
\epsfverbosetrue    % messages will show height and width
\centerline{\epsfbox{fig3.epsi}}
\bigskip

\head 3. Semiprojectivity\endhead

In this section we prove the main theorem of the paper, Theorem 3.12.

\proclaim{Theorem 3.1} Let $E$ be an infinite directed graph with
finite subgraphs $F_0\subseteq F_1\subseteq 
\ldots\subseteq E$ such that $\cup F_n = E$.  Let
$D:=\{u\in E^0\bigm|o\inv(u)\hbox{\ is\ infinite}\}$.  Suppose that 
 conditions (a) and (b) from Section 1 hold.
Then $\oh(E)$ is semiprojective.
\endproclaim

We will prove several lemmas before proving Theorem 3.1.  Lemmas 3.2, 
3.3, and 3.6 provide the main approximation arguments used in the 
proof of Theorem 3.1. Lemma 3.2 is a variation of 
a standard approximation result.  See, e.g., [B1, 2.23].

\proclaim{Lemma 3.2} There is $\delta>0$ and a map  $x\mapsto y\equiv 
y(x)$, defined on all elements of a \cstar-algebra for which $\Vert 
xx^*x-x\Vert<\delta$, having the following properties:
\item{(i)} $y$ is a partial isometry.
\item{(ii)} $y\in C^*(x)$.
\item{(iii)} $y-x$ lies in the ideal of $C^*(x)$ generated by 
$xx^*x-x$.
\item{(iv)} $\Vert y-x\Vert\to0$ as $\Vert xx^*x-x\Vert\to0$.
\item{(v)} If $p$ is a multiplier projection with $px=x$ (respectively, 
$px=0$, $xp=x$, or $xp=0$), then the same holds for $p$ and $y$.
\endproclaim

\demo{Proof}  Letting $f(t)=t^{-1/2}\Chi_{(1/2,\infty)}(t)$ for $t\ge0$, then 
$y(x)=xf(x^*x)$.
\qed
\enddemo

The next lemma uses Lemma 3.2 to `straighten' a family of approximate 
generators that extend a given set of generators from one subgraph to 
the next.  It is arranged so that we may use it on four separate 
occasions in the proof of Theorem 3.1.

\proclaim{Lemma 3.3}  Let $F\subseteq G$ be an inclusion of finite directed 
irreducible graphs, neither a cycle.  Let $D\subseteq F^0$, and 
suppose that
\smallskip
\item{(d1)} $o\inv(F^0\setminus D)\subseteq F^1$,
\item{(d2)} $o\inv(F^0)\cap t\inv(F^0)\subseteq F^1$.
\smallskip
Let $M$ be a directed set, and for $\mu\in M$ let $B_\mu$ be a 
\cstar-algebra with ideal $I_\mu$.  Suppose that
for each $\mu\in M$ we are given elements
$\bigl\{c\submu(e)\bigm|e\in G^1\bigr\}\subseteq B_\mu$
such that properties (d3) -- (d5) (below) hold:
\smallskip
\item{(d3)} $\bigl\{c\submu(e)\bigm|e\in F^1\bigr\}$ satisfy 
condition $(\oh)$.
\smallskip\noindent
For $u\in F^0$ let $\xi\submu(u)=c\submu(e)^*c\submu(e)$ for any 
$e\in F^1\cap t\inv(u)$.  Let 
$p\submu=\sum\bigl\{\xi\submu(u)\bigm| u\in F^0\bigr\}$, 
and for $u\in 
D$, let $q\submu(u)=\sum\bigl\{c\submu(e)c\submu(e)^*\bigm| e\in F^1\cap 
\,o\inv(u)\bigr\}$.
\smallskip
\item{(d4)} For $e\in G^1\setminus F^1$,  $c\submu(e)$ satisfies:
\smallskip
$$c\submu(e)=\cases
(1-p\submu)c\submu(e)\xi\submu\ofte,&
\text{if $o(e)\not\in D$ and $t(e)\in F^0$;}\cr
(1-p\submu)c\submu(e)(1-p\submu),&
\text{if $o(e)\not\in D$ and $t(e)\not\in F^0$;}\cr
\bigl(\xi\submu(u)-q\submu(u)\bigr)c\submu(e)(1-p\submu),&
\text{if $o(e)=u\in D$.}\cr\endcases$$
\smallskip
\item{(d5)} $\bigl\{c\submu(e)\bigm|e\in G^1\bigr\}$ satisfy 
condition $(\oh)$ asymptotically, in the sense that the relevant 
quantities belong to $I\submu$, and tend to 0 as $\mu\to\infty$.
\smallskip\noindent
Then there exists $\mu_0$ such that for all $\mu\ge\mu_0$ there are 
elements $\bigl\{a\submu(e)\bigm|e\in G^1\bigr\}\subseteq B_\mu$
for which
properties (d6) -- (d9) (below) hold:
\smallskip
\item{(d6)} $\bigl\{a\submu(e)\bigm|e\in G^1\bigr\}$ satisfy 
condition $(\oh)$,
\smallskip
\item{(d7)} Property (d4) above holds for 
$\bigl\{a\submu(e)\bigm|e\in G^1\setminus F^1\bigr\}$,
\smallskip
\item{(d8)} $a\submu(e)=c\submu(e)$ for $e\in F^1$,
\smallskip
\item{(d9)} For all $e\in G^1$, $a\submu(e)-c\submu(e)\in I\submu$ and 
$\Vert a\submu(e)-c\submu(e)\Vert\to0$ as $\mu\to\infty$.
\endproclaim

\demo{Proof} Let $G^1\setminus F^1\setminus 
o\inv(D)=\{e_1,e_2,\ldots,e_m\}$.  For all large enough $\mu$, 
$c\submu(e_1)$ satisfies the hypotheses of Lemma 3.2.  For such $\mu$,
let $b\submu(e_1)$ be the partial 
isometry obtained from $c\submu(e_1)$ as in Lemma 
3.2.  Inductively, for all large enough $\mu$, we can
let $b\submu(e_j)$ be the partial isometry obtained from 
$\bigl[1-\sum_{i<j}b\submu(e_i)b\submu(e_i)^*\bigr] 
c\submu(e_j)$ as in Lemma 3.2.  Then $\{b\submu(e_1),\ldots,
b\submu(e_m)\}$ 
are partial isometries whose final projections are pairwise orthogonal.
  For $u\in G^0\setminus 
F^0$ let 
$$\leqalignno{b\submu(u)&=\sum\bigl\{ b\submu(e)b\submu(e)^* \bigm| e\in 
o\inv(u)\bigr\}.\cr
\noalign{Note that}
b\submu(u)&\le 1-p\submu,&(d10)\cr}$$
by the first two lines of property (d4).

We now consider $(G^1\setminus F^1)\cap o\inv(D)$.  Fix $u\in D$. 
Let
$(G^1\setminus F^1)\cap o\inv(u) = \{f_1,\ldots,f_n\}$.  As before,
for all large enough $\mu$
let $b\submu(f_j)$ be the partial isometry obtained from 
$\bigl[1-\sum_{i<j}b\submu(f_i)b\submu(f_i)^*\bigr] 
c\submu(f_j)$ as in Lemma 3.2. Then $\{b\submu(f_1),
\ldots,\allowmathbreak b\submu(f_n)\}$ 
are partial isometries whose final projections are pairwise orthogonal.

Notice that property (d4) holds for $\bigl\{b\submu(e)\bigm|e\in 
G^1\setminus F^1\bigr\}$, by Lemma 3.2 and (d10).  Notice also that for 
all $e\in G^1\setminus F^1$, $b\submu(e)-c\submu(e)\in I\submu$ and 
$\Vert b\submu(e)-c\submu(e)\Vert\to0$ as $\mu\to\infty$, by (d5).

Now, for $e\in G^1\setminus F^1$, and all large enough $\mu$,
if $t(e)\not\in F^0$ then 
$b\submu(e)^*b\submu(e)-b\submu\bigl(t(e)\bigr)$
is a small element of $I\submu$ (by the asymptotic version of
condition (o6)).  Similarly, if $t(e)\in 
F^0$, then $b\submu(e)^*b\submu(e)-\xi\submu\bigl(t(e)\bigr)$ 
is a small element of $I\submu$.  It follows from Lemma 
3.2 that for all large enough $\mu$, there is a partial isometry 
$d\submu(e)\in B\submu$ such that
$$\eqalign{d\submu(e)^*d\submu(e)&
=\cases b\submu\bigl(t(e)\bigr),&
\text{if $t(e)\not\in F^0$,}\cr\cr
\xi\submu\bigl(t(e)\bigr),&\text{if $t(e)\in F^0$;}\cr\endcases
\cr\cr
d\submu(e)d\submu(e)^*&=b\submu(e)^*b\submu(e);\cr\cr
d\submu(e)- b\submu(e)^*b\submu(e)&\in I\submu.\cr\cr
\Vert d\submu(e)- b\submu(e)^*b\submu(e)\Vert&
\to0\hbox{\ as\ }\mu\to\infty.\cr\cr}$$
Set 
$$a\submu(e)=\cases b\submu(e)d\submu(e),&
\text{if $e\in G^1\setminus F^1$,}\cr\cr
c\submu(e),&\text{if $e\in F^1$. \qed}\cr\endcases$$
\enddemo

\definition{Notation 3.4} Let $B$ be a \cstar-algebra, $J$ an ideal in $B$, 
and $\Lambda$ a quasi-central approximate unit for $J$.  If $x=
x(h)$ and $y= y(h)$ are elements of $B$ depending on a variable 
$h\in\Lambda$, we will write $x\approx y$ to mean that $x-y\in J$ for 
all $h$, and $\lim_h\Vert x-y\Vert=0$.
\enddefinition

\remark{Remark 3.5}  Let $B$ be a \cstar-algebra, $J$ an ideal in $B$, 
and $\Lambda$ a quasi-central approximate unit for $J$.  For 
$h\in\Lambda$ let $k= k(h)=(1-h^2)^{1/2}$.  Then for any 
$x\in B$, $kx\approx xk$, and for any $y\in 
J$, $ky\approx0$.
\endremark

The next lemma is the key technique used in the proof of Theorem 3.1.  
It allows us to construct a partial isometry from a given partial 
isometry and a more suitable element which, however,
 only becomes a partial isometry after 
dividing by an ideal.

\proclaim{Lemma 3.6} Let $B$ be a \cstar-algebra and let $J$ be an ideal 
in $B$.  Let $a,\,b\in B$ be such that $a$ is a partial isometry
and $b+J$ is a partial isometry (in $B/J$).
Suppose that one of the following cases holds:
\smallskip\noindent
\item{(i)} $a^*b$, $a^*a-b^*b\in J$;
\smallskip\noindent
\item{(ii)} $a^*b$, $ab^* \in J$, and $a^*a+J\sim b^*b+J$ in $B/J$.
\smallskip\noindent
Let 
$r\in B$ be such that
$$r^*r-a^*a,\; rr^*-b^*b\ \in J,$$
where we let $r=a^*a$ in case (i).
Let $\Lambda$ be a quasi-central approximate unit for $J$.  Let $h$ 
and $k$ be as in Remark 3.5.  For $h\in \Lambda$ define $w= 
w(h)\in B$ by:
$$w=\cases ha+kb,&\text{in case (i);}\cr
h^2a+k^2b+hk(ar^*+br),&\text{in case (ii).}\cr\endcases$$
Then $ww^*w\approx w$.  Moreover we have the following asymptotic 
formulas for $w^*w$ and $ww^*$:
$$\eqalignno{w^*w&\approx\cases a^*a,&\text{in case (i);}\cr
h^2a^*a+k^2b^*b+hk(r+r^*),&\text{in case (ii);}\cr\endcases\cr
\noalign{\hbox{and}}
ww^*&\approx
h^2aa^*+k^2bb^*+hk(bra^*+ar^*b^*).\cr}$$
\endproclaim

\demo{Proof} We first consider case (i).  Using Remark 3.5 we find
$$\eqalign{w^*w&=a^*h^2a+b^*k^2b+2\real(a^*hkb)\cr
&\approx h^2a^*a+k^2b^*b+2\real(hka^*b)\cr
&\approx h^2a^*a+k^2a^*a+0\cr
&=a^*a.\cr}$$
This establishes the first asymptotic formula.  Moreover,
since $a^*a$ is a projection, we have $(w^*w)^2\approx w^*w$, and 
hence $(ww^*w-w)^*(ww^*w-w)=(w^*w)^3-2(w^*w)^2+w^*w\approx0$.  Thus 
$ww^*w\approx w$.
For the second asymptotic formula, we have
$$\eqalign{ww^*&=haa^*h+kbb^*k+2\real(hab^*k)\cr
&\approx h^2aa^*+k^2bb^*+hk(bra^*+ar^*b^*).\cr}$$

We now treat case (ii).  First note that, along with those explicitly 
mentioned in the statement of the Lemma,
the following elements (and 
their adjoints)
belong to $J$:  $ar$, $br^*$, $ra^*a-r$, $b^*br-r$,
 and $r^2$ (the last
because $r-rr^*r\in J$ and
 $(rr^*r)^2=ra^*ab^*br\in J$, 
so that $r^2=\bigl[rr^*r + (r-rr^*r)\bigr]^2\in J$).  
In particular, by Remark 3.5, 
if $x$ is any of these elements, then $kx\approx0$.  Using Remark 3.5 
repeatedly gives
$$\eqalignno{w^*w&=
\bigl[a^*h^2+b^*k^2+(ra^*+r^*b^*)hk\bigr] 
\bigl[h^2a+k^2b+hk(ar^*+br)\bigr]\cr
&\approx h^4a^*a + h^3ka^*ar^* + k^4b^*b + hk^3b^*br\cr
\noalign{\hskip1true in $+ h^3kra^*a + 
hk^3r^*b^*b + h^2k^2(ra^*ar^* + r^*b^*br)$}
&\approx h^4a^*a + h^3kr^* + k^4b^*b + hk^3r + h^3kr + hk^3r^* + 
h^2k^2(b^*b + a^*a)\cr
&= (h^4+h^2k^2)a^*a + (k^4+h^2k^2)b^*b + (h^3k+hk^3)(r+r^*)\cr
&= h^2a^*a + k^2b^*b + hk(r+r^*),\cr}$$
establishing the first asymptotic formula.  The second is proved in 
the same way, and we omit the details.  Finally, applying similar 
computations to the first asymptotic formula gives $(w^*w)^2\approx 
w^*w$, and hence, as in case (i), $ww^*w\approx w$.\qed
\enddemo

At this point we wish to briefly sketch the main idea of the proof of 
Theorem 3.1.  The first step is to lift an initial Toeplitz 
algebra, $\t\oh(F_2)$, modulo some ideal.  Then we extend that lift to 
a $*$-homomorphism on all of $\oh(E)$.  However, we will inductively 
extend the lifting from $F_0$ to $F_1$, etc.  In order to use Lemma 3.6 
for this
we have to arrange the initial extension so as to `miss' $F_1\setminus 
F_0$.  Lemma 3.9 below is our way of making this arrangement.

\proclaim{Lemma 3.7}  Let $F$ be an irreducible graph and not a 
cycle.  Let $D\subseteq F^0$ and $a\in F^0$.  Then $\t\oh(F,D^c)$ 
contains an infinite family of partial isometries 
with pairwise orthogonal final projections less than $a$, 
and initial projections equal to $a$.
\endproclaim

\demo{Proof} Choose two distinct cycles at $a$.  Let $s$ and $t$ 
denote the corresponding products of edges.  Then $s$, 
$t\in\t\oh(F,D^c)$ satisfy
$$\eqalign{s^*s=t^*t&=a,\cr
ss^*,\;tt^*&\le a,\cr
s^*t&=0.\cr}$$
It follows that $\bigl\{s^it\bigm|i=0,1,2,\ldots\bigr\}$ are partial 
isometries with initial projections all equal to $a$, and pairwise 
orthogonal final projections less than $a$. \qed
\enddemo

\proclaim{Lemma 3.8} Let $F$, $D$, and $a$ be as in Lemma 3.7.  Then 
$\t\oh(F,D^c)$ contains a family $\Omega$ of pairwise orthogonal 
subprojections of $a$ such that every vertex in $F^0$ is equivalent 
to infinitely many elements of $\Omega$.
\endproclaim

\demo{Proof} By Lemma 3.7 there is an infinite family of partial isometries in 
$\t\oh(F,D^c)$ with pairwise orthogonal final projections less than $a$, 
and initial projections equal to $a$.  Divide these into pairwise 
disjoint infinite subfamilies $\bigl\{\Delta_u\bigm|u\in F^0\bigr\}$.  For 
each $u\in F^0$ choose a path from $a$ to $u$; let $r_u$ denote
the corresponding product of edges.  Then $r_u^*r_u=u$ and 
$r_ur_u^*\le a$.  Let $\Omega=\cup_{u\in 
F^0}\bigl\{sr_ur_u^*s^*\bigm|s\in\Delta_u\bigr\}$. \qed
\enddemo

\proclaim{Lemma 3.9}  Let $E$ and $(F_n)_{n\ge0}$
satisfy conditions (a) and (b)
(and hence also (a6), (b3), and (b4), by Remarks 1.2).
For each $n\ge2$ choose $a_n\in F^0_n\setminus F^0_{n-1}$, and a map 
$e_n:D\to F^1_n\setminus F^1_{n-1}$ with $e_{n,u}\in o\inv(u)$ for 
$u\in D$.  (Note that $e_{n,u}$ exists by condition (a6) of Remarks 1.1.)
Then there is a $*$-homomorphism $\phi_n:\oh(E)\to\t\oh(F_n)$ such 
that
\smallskip
\halign{\quad\hfil#&\quad$#$\hfil&\quad for\quad$#$\hfil\cr
(c1)&\phi_n(u)\sim u&u\in E^0;\cr
(c2)&\phi_n(u)\le a_n&u\in E^0\setminus D;\cr
(c3)&\phi_n(u)=u&u\in D;\cr
(c4)&\phi_n(f)\phi_n(f)^*\le e_{n,u}e_{n,u}^*&f\in o\inv(u),\ u\in 
D.\cr}
\endproclaim

\demo{Proof} Fix $n\ge2$.  For this proof we will omit the subscript $n$ (so 
that $F=F_n$, $a=a_n$, $e_u=e_{n,u}$, etc.).  
By Lemma 3.8 there is a family $\Omega$ 
of pairwise orthogonal subprojections of $a$ in $\t\oh(F)$ such 
that each vertex in $F^0$ is equivalent (in 
$\t\oh(F)$) to infinitely many elements of
$\Omega$.  For each $u\in D$ choose a path 
in $F$ from $t(e_u)$ to $a$; let $\xi_u$ denote the corresponding 
product of edges.  Let 
$\Omega_u=\bigl\{e_u\xi_u\omega\xi_u^*e_u^*\bigm|\omega\in\Omega\bigr\}$.  
Then $\Omega_u\subseteq\t\oh(F)$ is a pairwise orthogonal family of 
subprojections of $e_ue_u^*$ such 
that each vertex in $F^0$ is equivalent (in 
$\t\oh(F)$) to infinitely many elements of
$\Omega_u$.

We first define $\phi$ on $E^0$, then on $E^1$.  For $u\in D$ we let 
$\phi(u)=u$.  For $u\not\in D$ choose $v_u\in F^0_0\setminus D$ as in 
(b2). (Let $v_u=u$ if $u\in F^0_0$.)  Thus there is a graph 
isomorphism $\gamma_u:o\inv(u)\to o\inv(v_u)$ such that for all $e\in 
o\inv(u)$, $\bigl[t(e)\bigr] = \bigl[t\bigl(\gamma_u(e)\bigr)\bigr]$ 
in $K_0\bigl(\oh(E)\bigr)$.  Choose $\phi:E^0\setminus D\to\Omega$ so 
that $\phi(u)\sim v_u$ in $\t\oh(F)$, and so that $\phi$ is 
one-to-one.  Let $\alpha_u\in\t\oh(F)$ with
$$\eqalign{\alpha_u^*\alpha_u&=\phi(u)\cr
\alpha_u\alpha_u^*&=v_u.\cr}$$
(For $u\in D$ define $\alpha_u=v_u=u$.)

We claim that $u\sim v_u$ in $\oh(E)$, for all $u\in E^0$.  This is 
immediate for $u\in D$.  For $u\in E^0\setminus D$ consider $e\in 
o\inv(u)$.  We have $\bigl[t(e)\bigr] = 
\bigl[t\bigl(\gamma_u(e)\bigr)\bigr]$.  By [C2], $t(e)\sim
t\bigl(\gamma_u(e)\bigr)$ in $\oh(E)$.  Hence
$$\eqalign{u
&=\sum\bigl\{ee^*\bigm|e\in o\inv(u)\bigr\}\cr
&\sim\sum\bigl\{\gamma_u(e)\gamma_u(e)^*\bigm|e\in o\inv(u)\bigr\}\cr
&=v_u,\cr}$$
establishing the claim.  Now if $e\in o\inv(u)$ for $u\not\in D$, then
$$\bigl[t\bigl(\gamma_u(e)\bigr)\bigr] = \bigl[t(e)\bigr]= 
\bigl[v_{t(e)}\bigr].$$
By (b4), $t\bigl(\gamma_u(e)\bigr)\sim v_{t(e)}$ in $\t\oh(F)$.  
Hence $t\bigl(\gamma_u(e)\bigr)\sim \phi\ofte$ in $\t\oh(F)$.  Let 
$\beta_e\in\t\oh(F)$ with
$$\eqalign{\beta_e^*\beta_e&=\phi\ofte\cr
\beta_e\beta_e^*&=t\bigl(\gamma_u(e)\bigr).\cr}$$

We define $\phi$ on $E^1\setminus o\inv(D)$ by
$$\displaylines{\quad\phi(e)=\alpha_u^*\gamma_u(e)\beta_e,
\hbox{\ for\ }e\in o\inv(u),\ u\not\in D.\cr}$$
We compute for $u\not\in D$:
$$\eqalign{\phi(e)^*\phi(e)&=\beta_e^*\gamma_u(e)^*\alpha_u 
\alpha_u^*\gamma_u(e)\beta_e\cr
&=\beta_e^*\gamma_u(e)^*v_u\gamma_u(e)\beta_e\cr
&=\beta_e^*t\bigl(\gamma_u(e)\bigr)\beta_e\cr
&=\phi\ofte,\cr\cr
\phi(e)\phi(e)^*&=\alpha_u^*\gamma_u(e)\beta_e\beta_e^* 
\gamma_u(e)^*\alpha_u\cr
&=\alpha_u^*\gamma_u(e)\gamma_u(e)^*\alpha_u,\cr\cr
\sum_{e\in o\inv(u)}\phi(e)\phi(e)^*&=
\alpha_u^*\left(\sum_{e\in o\inv(u)}\gamma_u(e)\gamma_u(e)^*\right)
\alpha_u\cr
&=\alpha_u^*v_u\alpha_u\cr
&=\phi(u).\cr}$$

Finally, let $u\in D$.  Choose a one-to-one map 
$\epsilon_u:o\inv(u)\to\Omega_u$ such that 
$\epsilon_u(e)\sim\phi\ofte$ in $\t\oh(F)$ for $e\in o\inv(u)$.  For 
$e\in o\inv(u)$ let $\phi(e)\in\t\oh(F)$ satisfy
$$\eqalign{\phi(e)^*\phi(e)&=\phi\ofte\cr
\phi(e)\phi(e)^*&=\epsilon_u(e).\cr}$$
It is now easy to see that $\phi:\oh(E)\to\t\oh(F)$ satisfies the 
requirements of the Lemma. \qed
\enddemo

Lemma 3.11 below is a necessary technical device for the main inductive 
step in the proof of Theorem 3.1.  We let $\u(A)$ denote the unitary 
group of a unital \cstar-algebra $A$.

\proclaim{Lemma 3.10}  Let $F$ be a finite irreducible graph which is not 
a cycle.  Let $D\subseteq F^0$, and let $\t\oh(F)$ denote 
$\t\oh(F,D^c)$.  Then for any vertex $a\in F^0$, the map
$$\u\bigl(a\cdot\t\oh(F)\cdot a\bigr)\to K_1\bigl(\t\oh(F)\bigr)$$ is 
onto.
\endproclaim

\demo{Proof} Let $x=(x(e))_{e\in F^1}\in\z^{|F^1|}$ represent an element of 
$K_1\bigl(\t\oh(F)\bigr)$.  Then $x(e)=0$ whenever $o(e)\in D$.  It 
follows that the element $U_0(x)\in M_{|x|}\bigl(\t\oh(F)\bigr)$
constructed in Lemma 2.1 is 
unitary, and $[U_0(x)]$ represents $x$ in 
$K_1\bigl(\t\oh(F)\bigr)$.

Now let $a\in F_0$.  
It follows from Lemma 3.7 that there are partial isometries $s_1$, 
$s_2$, $\ldots\in\t\oh(F)$ with initial projections equal to $a$ and 
pairwise orthogonal final projections less than $a$.
For each $u\in F^0$ choose a path 
from $a$ to $u$, and let $r_u$ denote the corresponding product of 
edges.  Then $r_u^*r_u=u$ and $r_ur_u^*\le a$.  Let $F^0 = 
\{u_1,u_2,\ldots,u_k\}$, and for $n\ge1$ set
$$t_n = \sum_{i=1}^k s_{(n-1)k+i}\,r_{u_i}.$$
Then $t_1$, $t_2$, $\ldots$ are isometries in $\t\oh(F)$ with 
orthogonal ranges dominated by $a$.  Now $t = (t_1,t_2,\ldots,t_k)$ 
is a partial isometry with
$$t^*t=1_{M_{|x|}\t\oh(F)}\hbox{\ and\ }tt^*\le a.$$
Then $tU_0(x)t^*+a-tt^*\in \u\bigl(a\cdot\t\oh(F)\cdot a\bigr)$ is 
equivalent to $U_0(x)$.\qed
\enddemo

\proclaim{Lemma 3.11}  Let $B$ be a \cstar-algebra, let $\l$ be a directed 
family of ideals of $B$ with closure $I$, and suppose that $B/I$ is 
simple and purely infinite.  Let $p\in B$ be a projection such that 
$\u(pBp)\to K_1(B/I)$ is onto.  Let $x\in B$ and $J_1\in\l$ be such that
$$x^*x-p,\;xx^*-p\in J_1.$$
Then there is $J_2\in\l$ with $J_1\subseteq J_2$, and there 
is $y\in B$ such that
$$y^*y=yy^*=p \qquad\hbox{and}\qquad x-y\in J_2.$$
\endproclaim

\demo{Proof} Replacing $B$ by $pBp$ and $\l$ by $p\l p$, we may assume that $B$ 
is unital and that $x+J_1$ is unitary.  By the assumption on $K$-theory 
we may choose $y_1\in\u(B)$ such that $[y_1+I]=[x+I]$ in $K_1(B/I)$.  
Then $[y_1^*x +I]=0$.  By [C2], $y_1^*x+I$ is in the 
connected component of the identy in $\u(B/I)$, and hence is a 
product of exponentials.
Since exponentials lift to 
$B$ there is $y_2\in\u(B)$ with $y_1^*x-y_2\in I$.  Then $x-y_1y_2\in 
I$, and so there is $J_2\in\l$, $J_1\subseteq J_2$, such 
that $\Vert x-y_1y_2+J_2\Vert<1$.  Now $y_2^*y_1^*x+J_2$ 
lies in the unit ball centered at the identity, and hence is an 
exponential.  Thus there is $y_3\in\u(B)$ with $y_2^*y_1^*x-y_3\in 
J_2$.  We may take $y=y_1y_2y_3$.\qed
\enddemo

\demo{Proof of Theorem 3.1}  Let $B$ be a \cstar-algebra,
 let $\l$ be a directed 
family of ideals of $B$ with closure $I$, and suppose that $B/I$ is 
isomorphic to $\oh(E)$ (we will identify $B/I$ with $\oh(E)$).  
For any $J\in\l$, we
will use $\pi:B/J\to B/I$ to denote the quotient map (the context 
should make it clear what the domain of $\pi$ is).  If $J^\prime\subseteq 
J$ are two elements of $\l$, we will write 
$\pi_J:B/J^\prime\to B/J$, or $\pi_n:B/J_i\to B/J_n$ if $J_i\subseteq 
J_n$, and if $x'\in J'$, $x\in J$, we will write $x'x$ for $\pi_J(x')x$. 
We also recall the $*$-homomorphisms 
$\phi_n:\oh(E)\to\t\oh(F_n)$ provided by Lemma 3.9.

Note that Lemma 3.3 is applicable in the case where the smaller graph 
is empty.  Thus 
we can find $J_0\in\l$ and a $*$-homomorphism 
$\theta_0:\t\oh(F_2)\to B/J_0$ such that
\halign{$#$\hfil&\qquad$#$\hfil&\quad#\hfil\cr
\noalign{\smallskip}
(2_0)&\pi\circ\theta_0=\hbox{id}_{\t\oh(F_2)}.\cr
\noalign{\smallskip}
\noalign{\noindent
For $u\in F^0_0\setminus D$, recall from the 
proof of Lemma 
3.9 the element $\alpha_{2,u}\in\t\oh(F_2)$ with 
$\alpha_{2,u}^*\alpha_{2,u} = \phi_2(u)$ and 
$\alpha_{2,u}\alpha_{2,u}^* = v_u = u$.  Define $\psi_0:\oh(E)\to 
B/J_0$ by:}
\noalign{\smallskip}
(3_0)&\psi_0\restrictedto{\t\oh(F_0)} = 
\theta_0\restrictedto{\t\oh(F_0)};\cr
\noalign{\smallskip}
(4_{00})&\psi_0(e)=\theta_0\bigl(\phi_2(e)\bigr),&if $e\in E^1\setminus 
F_0^1$ and $t(e)\not\in F_0^0$;\cr\cr
&\psi_0(e)=\theta_0\bigl(\phi_2(e)\alpha_{2,t(e)}^*\bigr),&if $e\in 
E^1\setminus F_0^1$ and $t(e)\in F_0^0$.\cr
\noalign{\smallskip}
\noalign{\noindent
It is easy to check that $\psi_0$ is a $*$-homomorphism. We 
remark that $(4_{00})$ implies the following:}
\noalign{\smallskip}
(4_0)&\psi_0(u)=\theta_0\bigl(\phi_2(u)\bigr),&if $u\in 
E^0\setminus F^0_0$.\cr\cr
&\psi_0(ee^*)=\theta_0\bigl(\phi_2(ee^*)\bigr),&if $e\in 
E^1\setminus F^1_0$ and $o(e)\in D$.\cr}
\medskip
Suppose inductively that we have constructed for $j<n$ the following:
$$\displaylines{
\quad J_j\in\l\hfill\cr
\quad \theta_j:\t\oh(F\jplpl)\to B/J_j,\hbox{\ a\ 
}*\hbox{-homomorphism} \hfill\cr
\quad\psi_j:\oh(E)\to B/J_0,\hbox{\ a\ 
}*\hbox{-homomorphism} \hfill\cr}$$
such that
\halign{$#$\hfil&\qquad$#$\hfil&\quad#\hfil\cr
\noalign{\smallskip}
(1_j)&J\jmi\subseteq J_j\cr
\noalign{\smallskip}
(2_j)&\pi\circ\theta_j=\hbox{id}_{\t\oh(F\jplpl)}\cr
\noalign{\smallskip}
(3_j)&\pi_j\circ\psi_j\restrictedto{\t\oh(F_j)}=\theta_j 
\restrictedto{\t\oh(F_j)}\cr
\noalign{\smallskip}
(4_j)&\pi_j\circ\psi_j(u)=\theta_j\bigl(\phi\jplpl(u)\bigr),& 
if $u\in E^0\setminus F^0_j$\cr
\noalign{\smallskip}
&\pi_j\circ\psi_j(ee^*)=\theta_j\bigl(\phi\jplpl(ee^*)\bigr),&
if $e\in E^1\setminus F^1_j$ and $o(e)\in D$\cr
\noalign{\smallskip}
(5_j)&\theta_j\restrictedto{\t\oh(F\jpl)}=\pi_j\circ\theta\jmi\cr
\noalign{\smallskip}
(6_j)&\psi_j\restrictedto{\t\oh(F\jmi)}=\psi\jmi 
\restrictedto{\t\oh(F\jmi)}.\cr}

\smallskip\noindent
We will construct $J_n$, $\theta_n$, $\psi_n$ satisfying the above 
for $j=n$.  Then $\psi=\lim_n\psi_n$ will establish the theorem.
First we will define $\theta_n$.  For this we need to consider edges 
in $F^1\nplpl\setminus F^1\npl$.

For $e\in F^1\nplpl\setminus F^1\npl$ choose $\sigma(e)\in B/J\nmi$ 
such that
\smallskip
\halign{$#$\hfil&\qquad$#$\hfil&\quad#\hfil\cr
\noalign{\smallskip}
(7\nmi)&\pi\bigl(\sigma(e)\bigr)=e.\cr}
\smallskip\noindent
Let 
$$\eqalign{
p&=\sum\bigl\{\theta\nmi(u)\bigm|u\in F^0\npl\bigr\},\cr
q(u)&=\sum\bigl\{\theta\nmi(e)\theta\nmi(e)^*\bigm| e\in F^1\npl\cap 
o\inv(u)\bigr\},\hbox{\ for\ }u\in D.\cr}$$
We will apply Lemma 3.3 with
$$\eqalign{F&=F\npl\cr
G&=F\nplpl\cr
M&=\{J\in\l\mid J\supseteq J\nmi\}\cr
B_J&=B/J\cr
I_J&=I/J\cr
c_J(e)&=\pi_J\bigl(c(e)\bigr),\cr}$$
where
$$\eqalignno{
 c(e)&=\theta\nmi(e),\hbox{\ if\ }e\in F^1\npl,\cr\cr
\noalign{\hbox{and\ for\ $e\in F^1\nplpl\setminus F^1\npl$,}}
 \cr
 c(e)&=\cases 
  (1-p)\sigma(e)\theta\nmi\ofte,&
  \text{if $o(e)\not\in D$ and $t(e)\in F^0\npl$,}\cr\cr
  (1-p)\sigma(e)(1-p),&
  \text{if $o(e)\not\in D$ and $t(e)\not\in F^0\npl$,}\cr\cr
  \bigl(\theta\nmi(u)-q(u)\bigr)\sigma(e)(1-p),&
  \text{if $o(e)=u\in D$.}\cr\endcases\cr}$$
(Note that $p_J$ and $q_J(u)$ in Lemma 3.3 are given by $\pi_J(p)$ 
and $\pi_J\bigl(q(u)\bigr)$.)
Then properties (d1) -- (d4) of Lemma 3.3 are immediate, and (d5) follows 
from $(2\nmi)$ and $(7\nmi)$ above.  By Lemma 3.3 there are $J'_n\in\l$ 
and a $*$-homomorphism
$\theta'_n:\t\oh(F\nplpl)\to B/J'_n$ satisfying $(1_n)$, $(2_n)$, 
and $(5_n)$.

The construction of $\psi_n$ proceeds in two stages.  First we extend
$\psi\nmi$ from $F\nmi$ to $F_n$ satisfying $(3_n)$;
then we extend it from $F_n$ to 
$E$ satisfying $(4_n)$.
Each stage will require another application of Lemma 3.3.

For $u\in 
F_n^0\setminus F\nmi^0$ choose $r_0(u)\in B/J'_n$ with
$$\eqalign{r_0(u)^*r_0(u) &=\pi'_n\circ\psi\nmi(u)\cr
r_0(u)r_0(u)^*&=\theta'_n(u).\cr}\leqno(*)$$
(Such elements exist because 
$$\eqalign{\pi'_n\circ\psi\nmi(u)&=\pi'_n\circ\pi\nmi\circ\psi\nmi(u)\cr
&=\pi'_n\circ\theta\nmi\bigl(\phi\npl(u)\bigr), 
\hbox{\ \ by\ }(4\nmi),\cr
&=\theta'_n\bigl(\phi\npl(u)\bigr),\hbox{\ \ by\ }(5_n),\cr
&\sim\theta'_n(u),\hbox{\ \ by\ Lemma\ 3.9, (c1)}.)\cr}$$
For $u\in F^0\nmi$ we set $r_0(u)=\pi'_n\circ\psi\nmi(u)$.  For $u\in 
F_n^0\setminus F^0\nmi$ let
$$r(u)=\sum\bigl\{\theta'_n(e)r_0\ofte\psi\nmi(e)^*\bigm| e\in 
F_n^1\cap o\inv(u)\bigl\},$$
and set $r(u)=r_0(u)$ for $u\in F^0\nmi$.  Then $(*)$ holds for 
$r(u)$ in place of $r_0(u)$.  Let $u\in F_n^0\setminus F^0\nmi$.  By $(*)$ 
for $r_0(u)$ and for $r(u)$, we have that $r_0(u)^*r(u)$ is a unitary 
on (the range of) $\pi'_n\circ\psi\nmi(u)$.  
By condition (b3) and Lemmas 3.10 and 3.11,
there are 
$y(u)\in B/J_0$ and $J(u)\in\l$ with $J'_n\subseteq J(u)$ such that
$$\eqalign{y(u)^*y(u)&=y(u)y(u)^*=\psi\nmi(u)\cr
\pi'_n\bigl(y(u)\bigr)&\equiv r_0(u)^*r(u) \pmod{J(u)}.\cr}$$
For $u\in F\nmi^0$ we let $y(u)=\psi\nmi(u)$.  Now for any $e\in 
F_n^1$ we set
$$\psi\nmi'(e)=\psi\nmi(e)y\ofte.$$
Then $\psi'\nmi(e)$ is a partial isometry having the same initial and 
final projections as $\psi\nmi(e)$.  
Hence $\bigl\{\psi'\nmi(e)\bigm| e\in F^1_n\bigr\}$ define a 
$*$-homomorphism of $\t\oh(F_n)$.
We note that 
$$\psi'\nmi\restrictedto{\t\oh(F\nmi)}=
\psi\nmi\restrictedto{\t\oh(F\nmi)}.$$

Choose $J_n\in\l$ with $J(u)\subseteq J_n$ for all $u\in 
F_n^0\setminus F^0\nmi$. We now define $\theta_n:\t\oh(F\nplpl)\to B/J_n$ by
$$\theta_n=\pi_n\circ\theta'_n.$$

Now let $u\in F_n^0\setminus F\nmi^0$ and $e\in o\inv(u)$.  We 
compute:
$$\eqalign{
\pi_n\bigl(r(u)\psi'\nmi(e)\psi'\nmi(e)^*\bigr)&=
  \pi_n\bigl( r(u)\psi\nmi(e)\psi\nmi(e)^*\bigr)\cr
 &=\pi_n\bigl(\theta_n'(e)r_0\ofte\psi\nmi(e)^*\bigr),
 \hbox{\ \ by\ definition\ of\ 
 }r(u),\cr
 &=\theta_n(e)r_0\ofte y\ofte\psi'\nmi(e)^*\cr
 &=\theta_n(e)r\ofte\psi'\nmi(e)^*.\cr}$$
Since
$$r(u)=\sum_{e\in o\inv(u)}r(u)\psi'\nmi(e)\psi'\nmi(e)^*,$$
then
$$\pi_n\bigl(r(u)\bigr)=\sum_{e\in o\inv(u)}
\theta_n(e)r\ofte\psi'\nmi(e)^*. \leqno(**)$$
For $u\in F_n^0$ and $e\in F_n^1$ choose arbitrarily $\widetilde 
r(u)$, $\widetilde\theta_n(u)$, $\widetilde\theta_n(e)\in B/J_0$ 
lifting $r(u)$, $\theta_n(u)$, and $\theta_n(e)$.  Let $\Lambda$ be a 
quasi-central approximate unit for $J_n$, and let $h$, $k$ be as in 
Remark 3.5.  For $e\in F_n^1\setminus F^1\nmi$ let $w(e)$ be defined as 
in Lemma 3.6, with $a=\psi'\nmi(e)$, $b=\widetilde\theta_n(e)$, and 
$r=\widetilde r\ofte$.  We will apply Lemma 3.3 with
$$\eqalignno{F&=F\nmi\cr
G&=F_n\cr
M&=\Lambda\cr
B_h&=B/J_0,\ h\in\Lambda\cr
I_h&=J_n/J_0,\ h\in\Lambda\cr
c_h(e)&=\psi\nmi(e),\hbox{\ if\ }e\in F^1\nmi,\cr\cr
\noalign{\hbox{and\ for\ $e\in F_n^1\setminus F^1\nmi$,}}
\cr
c_h(e)&=\cases (1-p_h)w(e)\psi\nmi\ofte,&
\text{if $o(e)\not\in D$ and $t(e)\in F^0\nmi$,}\cr\cr
(1-p_h)w(e)(1-p_h),&
\text{if $o(e)\not\in D$ and $t(e)\not\in F^0\nmi$,}\cr\cr
\bigl(\psi\nmi(u)-q_h(u)\bigr)w(e)(1-p_h),&
\text{if $o(e)=u\in D$,}\cr\endcases\cr}$$
where, as in the statement of Lemma 3.3, 
$p_h=\sum\bigl\{\psi\nmi(u)\bigm|u\in F^0\nmi\bigr\}$, and for $u\in 
D$, $q_h(u)=\sum\bigl\{\psi\nmi(e)\psi\nmi(e)^*\bigm|e\in F^1\nmi\cap 
o\inv(u)\bigr\}$.
Then (d1) -- (d4) of Lemma 3.3 hold by definition. To establish (d5) of 
Lemma 3.3, first note that by Remark 3.5, if $e\in F_n^1\setminus F^1\nmi$ 
then $c_h(e)\approx w(e)$.  

We will check that condition $(\oh)$ holds 
asymptotically.  It follows
from Lemma 3.6 that $w(e)w(e)^*w(e)\approx w(e)$, so 
$c_h(e)$ is asymptotically a partial isometry.  For $u\in 
F_n^0\setminus F^0\nmi$ and $e\in t\inv(u)$, Lemma 3.6 gives
$$\eqalign{w(e)^*w(e)&\approx h^2\psi\nmi(u) + 
k^2\widetilde\theta_n(e)^*\widetilde\theta_n(e) + 2hk\real\bigr( 
\widetilde r(u)\bigr)\cr
&\approx h^2\psi\nmi(u)+k^2\widetilde\theta_n(u) + 2hk\real\bigr( 
\widetilde r(u)\bigr),\cr}$$
which is independent of $e\in t\inv(u)$.  Define
$$c_h(u)=h^2\psi\nmi(u)+k^2\widetilde\theta_n(u) + 2hk\real\bigr( 
\widetilde r(u)\bigr).$$
Then $c_h(u)$ is asymptotically a projection.  This establishes (o1).

For (o2), let $u,v\in F_n^0$ be distinct.  If $u,v\in F^0\nmi$, then
$$\eqalignno{c_h(u)c_h(v)&=\psi\nmi(u)\psi\nmi(v)=0.\cr\cr
\noalign{\hbox{If\ $u\in F^0\nmi$\ and\ $v\not\in F^0\nmi$,}}
\cr
c_h(u)c_h(v)&=\psi\nmi(u)\,
\bigl[h^2\psi\nmi(v)+k^2\widetilde\theta_n(v) + 2hk\real\bigr( 
\widetilde r(v)\bigr)\bigr]\cr
&\approx0\cr\cr
\noalign{\hbox{by\ Remark\ F,\ since}}
\cr
\pi_n\circ\psi\nmi(u)&=\pi_n\circ\pi\nmi\circ\psi\nmi(u)\cr
&=\pi_n\circ\theta\nmi\bigl(\phi\npl(u)\bigr),\hbox{\ by\ }(5\nmi),\cr
&=\theta_n\bigl(\phi\npl(u)\bigr)\cr}$$
is orthogonal to $\theta_n(v)$.  The case where $u$, $v\not\in 
F^0\nmi$ is similar.

Condition (o3) follows from the definition of $c_h(u)$, and Lemma 3.6.

For (o4) let $u\in D$ and $e\in F^1_n\cap o\inv(u)$.  If $e\in 
F^1\nmi$, the condition is clear.  If $e\not\in F^1\nmi$, then
$$\eqalign{c_h(u)c_h(e)&\approx \psi\nmi(u)\bigl[h^2\psi'\nmi(e) + 
k^2 \widetilde\theta_n(e) + 2hk\real\bigl(\widetilde\theta_n(e) 
\widetilde r\ofte \psi'\nmi(e)^*\bigr)\bigr]\cr
&\approx c_h(e),\cr}$$
since $\psi\nmi(u)\psi'\nmi(e)=\psi'\nmi(e)$ and $\pi_n\bigl( 
\psi\nmi(u)\bigr)\theta_n(u)\theta_n(e)=\theta_n(e)$.

For (o5) let $u\in D$ and $e$, $f\in o\inv(u)\cap F_n^1$ with 
$e\not=f$.  If $e$, $f\in F^1\nmi$, the condition is clear.  If $e\in 
F^1\nmi$ and $f\not\in F^1\nmi$, then
$$c_h(e)^*c_h(f)=\bigl[q_h(u)c_h(e)\bigr]^*\bigl[ 
\bigl(\psi\nmi(u)-q_h(u)\bigr)c_h(f)\bigr]=0.$$
If $e$, $f\not\in F^1\nmi$ then $c_h(e)^*c_h(f)\approx0$, since
$$\eqalign{\pi_n\bigl(\psi'\nmi(e)^*\bigr)\,\theta_n(f) 
&=\pi_n\bigl(\psi\nmi(e)y\ofte\bigr)^*\,\theta_n(f)\cr
&=\pi_n\bigl(y\ofte\bigr)^*\, 
\pi_n\circ\theta\nmi\bigl(\phi\npl(e)\bigr)^*\, \theta_n(f)\cr
&=\pi_n\bigl(y\ofte\bigr)^*\,\theta_n\bigl(\phi\npl(e)\bigr)^*\,\theta_n(f)\cr
&=0.\cr}$$

For (o6) let $u\not\in D$.  If $u\in F^0\nmi$, the condition is 
clear.  If $u\in F^1_n\setminus F^1\nmi$ then $o\inv(u)\subseteq 
F_n^1\setminus F^1\nmi$.  We compute (notice the use of 
condition $(**)$):
$$\eqalign{
\sum_{e\in o\inv(u)} c_h(e)c_h(e)^*
&\approx\sum_{e\in o\inv(u)}w(e)w(e)^*\cr
&\approx\sum_{e\in o\inv(u)} h^2\psi'\nmi(e)\psi'\nmi(e)^* + 
k^2\widetilde\theta_n(e)\widetilde\theta_n(e)^*\cr
&\hbox{\hskip1.5true in}+ 2hk\real\bigl[ 
\widetilde\theta_n(e)\widetilde r\ofte \psi'\nmi(e)^*\bigr]\cr\cr
&\approx h^2\psi'\nmi(u) + k^2\widetilde\theta_n(u) + 2hk\real\bigl[ 
\widetilde r(u)\bigr],\qquad\hbox{\ by\ }(**),\cr
&=c_h(u).\cr}$$

It now follows from Lemma 3.3 that there are $\bigl\{\psi_n(e)\bigm| 
e\in F^1_n\bigr\}\subseteq B/J_0$ such that $(3_n)$ and $(6_n)$ hold, 
and $\psi_n:\t\oh(F_n)\to B/J_0$ is a $*$-homomorphism.  
Moreover, we may assume that
$$\psi_n(u)\sim\psi\nmi(u),\quad u\in F^0_n\setminus F^0\nmi.$$
The reason is that Lemma 3.6 (i) applies to the element
$$\lambda_u=h\psi\nmi(u)+k\widetilde r(u),$$
and so
$$\eqalign{\lambda_u^*\lambda_u&\approx\psi\nmi(u)\cr
\lambda_u\lambda_u^*&\approx c_h(u).\cr}$$
Thus in the previous application of Lemma 3.3, we may choose $h$ large 
enough in $\Lambda$ so that Lemma 3.2 applies to 
$\psi_n(u)\lambda_n\psi\nmi(u)$, giving the desired equivalence.

It remains 
to define $\bigl\{\psi_n(e)\bigm| e\in E^1\setminus F^1_n\bigr\}$ so 
that $\bigl\{\psi_n(e)\bigm|e\in E^1\bigr\}$ satisfy $(\oh)$ and so 
that $(4_n)$ holds.
We  will extend the definition of $r(u)$ to all $u\in E^0$,
and will define additional elements $s(e)$ for $e\in E^1\setminus 
F^1_n$ with $o(e)\in D$.
Recall from the proof of
Lemma 3.9 the elements $\alpha_{n,u}\in\t\oh(F_n)$, for $u\in 
E^0\setminus F^0_n$, with
$$\eqalign{\alpha_{n,u}^*\alpha_{n,u}&=\phi_n(u)\cr
\alpha_{n,u}\alpha_{n,u}^*&=v_u.\cr}$$
We let
$$\eqalign{r(u)&=\theta_n(\alpha_{n+2,u}^*\alpha_{n+1,u}), \quad 
u\in E^0\setminus F^0_n,\cr
s(e)&=\theta_n\bigl(\phi\nplpl(e)\alpha_{n+2,t(e)}^* 
\alpha_{n+1,t(e)} \phi\npl(e)^*\bigr),\quad e\in E^1\setminus F_n^1 
\hbox{\ and\ }o(e)\in D.\cr}$$
Then we have
$$\eqalign{r(u)^*r(u)&=\theta_n\bigl(\phi\npl(u)\bigr)=
 \pi_n\bigl(\psi\nmi(u)\bigr)\cr
r(u)r(u)^*&=\theta_n\bigl(\phi\nplpl(u)\bigr)\cr
s(e)^*s(e)&=\theta_n\bigl(\phi\npl(ee^*)\bigr)=
 \pi_n\bigl(\psi\nmi(ee^*)\bigr)\cr
s(e)s(e)^*&=\theta_n\bigl(\phi\nplpl(ee^*)\bigr).\cr}$$
Let $\widetilde r(u)$, $\widetilde s(e)\in B/J_0$ 
with $\pi_n\bigl( \widetilde 
r(u)\bigr) = r(u)$, $\pi_n\bigl(\widetilde s(u)\bigr)=s(u)$.  

Now we proceed with the definition of $\psi_n$.  Suppose inductively 
that for some $i>n$ we have extended 
$\psi_n\restrictedto{\t\oh(F_{i-2})}$ 
to a $*$-homomorphism 
$\psi_n:\t\oh(F\imi) \to B/J_0$ satisfying $(4_n)$ on 
$F\imi\setminus F_n$, and such that for $u\in F^0\imi\setminus 
F^0_{i-2}$ there exists $\lambda_u\in B/J_0$ with
$$\eqalign{\lambda_u^*\lambda_u&=\psi\nmi(u)\cr
\lambda_u\lambda_u^*&=\psi_n(u).\cr}$$
Let $\hi\in\Lambda$ and $\ki=(1-\hi^2)^{-1/2}$ be as in Remark 3.5.
We define
$$w_u=\hi\psi\nmi(u) + \ki\widetilde r(u),\quad 
u\in F_i^0\setminus F\imi^0,$$
and
$$w_e=\hi\psi\nmi(ee^*) + \ki\widetilde s(e),\quad 
e\in(F_i^1\setminus F^1\imi)\cap o\inv(D).$$
It follows from Lemma 3.6 (i) that
$$w_u^*w_u\approx\psi\nmi(u),\quad u\in F^0_i\setminus F^0\imi,$$
and
$$w^*_ew_e\approx\psi\nmi(ee^*),\quad
e\in(F_i^1\setminus F^1\imi)\cap o\inv(D).$$
For $e\in F^1_i\setminus F^1\imi$ we set
$$\eta(e)=\cases
w_{o(e)}\psi\nmi(e)\lambda_{t(e)}^*,&
\text{if $o(e)\not\in D$ and $t(e)\in F^0\imi$,}\cr
w_{o(e)}\psi\nmi(e)w^*_{t(e)},&
\text{if $o(e)\not\in D$ and $t(e)\not\in F^0\imi$,}\cr
w_e\psi\nmi(e)w_{t(e)}^*,&
\text{if $o(e)\in D$.}\cr\endcases$$
We will apply Lemma 3.3 with
$$\eqalignno{
F&=F\imi\cr
G&=F_i\cr
M&=\Lambda\cr
B_\hi&=B/J_0,\ \hi\in\Lambda\cr
I_\hi&=J_n/J_0,\ \hi\in\Lambda\cr
c_\hi(e)&=\psi_n(e),\hbox{\ if\ }e\in F^1\imi,\cr\cr
\noalign{\hbox{and\ for\ $e\in F^1_i\setminus F^1\imi$,}}
\cr
c_\hi(e)&=\cases
 (1-p_\hi)\eta(e)\psi_n\ofte,&
 \text{if $o(e)\not\in D$ and $t(e)\in F^0\imi$,}\cr\cr
 (1-p_\hi)\eta(e)(1-p_\hi),&
 \text{if $o(e)\not\in D$ and $t(e)\not\in F^0\imi$,}\cr\cr
 \bigl(\psi_n(u)-q_\hi(u)\bigr)\eta(e)(1-p_\hi),&
 \text{if $o(e)=u\in D$,}\cr\endcases\cr}$$
(where
$$\eqalign{p_\hi&=\sum\bigl\{\psi_n(u)\bigm| u\in F^0\nmi\bigr\}, 
\hbox{\ and}\cr
q_\hi(u)&=\sum\bigl\{\psi_n(e)\psi_n(e)^*\bigm| e\in F^1\nmi\cap 
o\inv(u)\bigr\},\quad u\in D,\cr}$$
as in Lemma 3.3.)
Then (d1) -- (d4) of Lemma 3.3 hold by definition.  Moreover, for $e\in 
F^1_i\setminus F^1\imi$ we have
$$c_\hi(e)\approx \eta(e),$$
as can be seen from the facts:
$$\eqalign{\pi_n(w_u)&=r(u),\quad u\in F^0_i\setminus F^0\imi,\cr
\noalign{\smallskip}
\pi_n(w_e)&=s(e),\quad e\in(F^1_i\setminus F^1\imi)\cap 
o\inv(u),\quad u\in D.\cr}$$

We will now check the asymptotic condition $(\oh)$ for 
$\bigl\{c_\hi(e)\bigm|i\in F^1_i\setminus F^1\imi\bigr\}$.  For (o1) 
we note that if $u\in F^0_i\setminus F^0\imi$ and $t(e)=u$, then
$$\eqalign{c_\hi(e)^*c_\hi(e)
&\approx \eta(e)^*\eta(e)\cr
&=w_{t(e)}\psi\nmi(e)^*w_{o(e)}^*
 w_{o(e)}\psi\nmi(e)w_{t(e)}^*\cr
&\approx w_{t(e)}\psi\nmi(e)^*\psi\nmi\ofoe
 \psi\nmi(e)w_{t(e)}^*\cr
&=w_{t(e)}\psi\nmi\ofte w_{t(e)}^*\cr
&\approx w_{t(e)}w_{t(e)}^*\cr}$$
is independent of $e\in t\inv(u)$.  Thus we may set
$$c_\hi(u)=w_u w_u^*.$$

For (o2) let $u$, $v\in F^0_i$ be distinct.  If $u$, $v\in F^0\imi$ 
then $c_\hi(u)c_\hi(v)=0$ by the inductive hypothesis.  If $u\in 
F^0\imi$ and $v\not\in F^0\imi$, then
$$\eqalign{c_\hi(u)c_\hi(v)
&\approx c_\hi(u)c_\hi(e)^*c_\hi(e),\hbox{\ for\ any\ } e\in 
t\inv(v),\cr
&=0,\cr}$$
by definition of $c_\hi(e)$.  If $u$, $v\not\in F^0\imi$, then we have
$$\eqalign{w_u^*w_v&\approx \hi^2\psi\nmi(u)\psi\nmi(v) + 
\hi\ki\bigl(\psi\nmi(u)\widetilde r(v) + \widetilde 
r(u)^*\psi\nmi(v)\bigr) + \ki^2\widetilde r(u)^*\widetilde r(v)\cr
&\approx0,\cr}$$
since $\psi\nmi(u)\widetilde r(v)$, $\psi\nmi(v)\widetilde r(u)$, 
and $\widetilde r(u)^*\widetilde r(v)\in J_n$.  Hence
$$c_\hi(u)c_\hi(v)=w_u w_u^*w_v w_v^* \approx0.$$

Condition (o3) is immediate from the definition of $c_\hi(u)$.

For (o4), let $u\in D$ and $e\in F^1_i\cap o\inv(u)$.  If $e\in 
F^1\imi$ the condition follows from the inductive hypothesis.  Suppose 
that $e\not\in F^1\imi$.  Then
$$\eqalignno{\psi_n(u)w_e&=\psi_n(u)\bigl(\hi\psi\nmi(ee^*) + 
\ki\widetilde s(e)\bigr)\cr
&\approx \hi\psi\nmi(ee^*) + \ki\widetilde s(e),\cr
\noalign{\hbox{since\ $\phi_\ell(ee^*)\le e_{\ell,u}e_{\ell,u}^* \le 
u$\ for\ all\ $\ell$,}}
&=w_e.\cr}$$
Since $u\in D$, $c_\hi(u)=\psi_n(u)$.  Then
$$\eqalign{c_\hi(u)c_\hi(e)
&\approx \psi_n(u)\eta(e)\cr
&\approx \psi_n(u)w_e\psi\nmi(e)w_{t(e)}^*\cr
&\approx w_e\psi\nmi(e)w_{t(e)}^*\cr
&=\eta(e)\cr
&\approx c_\hi(e).\cr}$$

For (o5), let $u\in D$ and $e$, $f\in F^1_i\cap o\inv(u)$ with 
$e\not=f$.  If $e$, $f\in F^1\imi$, the condition follows from the 
inductive hypothesis.  If $e\in F^1\imi$ and $f\not\in F^1\imi$, the 
condition follows from the definition of $c_\hi(f)$.  If $e$, 
$f\not\in F^1\imi$, then
$$\eqalignno{c_\hi(e)^*c_\hi(f)&\approx \eta(e)^*\eta(f)\cr
&= w_{t(e)}\psi\nmi(e)^*w_e^*w_fpsi\nmi(f)w_{t(f)}^*\cr
&\approx0,\cr
\noalign{\hbox{since}}
w_e^*w_f&\approx\bigl(\hi\psi\nmi(ee^*) + \ki\widetilde s(e)\bigr)^* 
\bigl(\hi\psi\nmi(ff^*) + \ki\widetilde s(f)\bigr)\cr
&\approx\hi\ki\bigl(\psi\nmi(ee^*)\widetilde s(f) + \widetilde 
s(e)^* \psi\nmi(ff^*)\bigr) + \ki^2\widetilde s(e)^*\widetilde s(f)\cr
&\approx0,\cr}$$
by Remark 3.5.

For (o6), let $u\in F^0_i\setminus D$.  If $u\in F^0\imi$ the 
condition follows from the inductive hypothesis.  If $u\not\in 
F^0\imi$, then
$$\eqalign{\sum\bigl\{c_\hi(e)&c_\hi(e)^*\bigm|e\in o\inv(u)\bigr\}
\approx\sum\bigl\{\eta(e)\eta(e)^*\bigm|e\in o\inv(u)\bigr\}\cr
&=\sum\bigl\{
 w_u\psi\nmi(e)w_{t(e)}^*w_{t(e)}
\psi\nmi(e)^*w_u^*\bigm|
 e\in o\inv(u)\setminus t\inv(F^0\imi)\bigr\}\cr
&\quad+ \sum\bigl\{w_u\psi\nmi(e)\lambda_{t(e)}\lambda_{t(e)}^* 
\psi\nmi(e)^*w_u^*\bigm|e\in o\inv(u)\cap t\inv(F^0\imi)\bigr\}\cr
&\approx\sum\bigl\{ w_u\psi\nmi(e)
\psi\nmi\ofte\psi\nmi(e)^*w_u^*\bigm|e\in o\inv(u)\bigr\}\cr
&=w_u\psi\nmi(u)w_u^*\cr
&\approx w_uw_u^*\cr
&=c_\hi(u).\cr}$$
It now follows from Lemma 3.3 that there are $\bigl\{\psi_n(e) \bigm| 
e\in F^1_i\setminus F^1\imi\bigr\}$ such 
that $\bigl\{\psi_n(e)\bigm|e\in F^1_i\bigr\}$ satisfy $(\oh)$.
We now check $(4_n)$ on $F_i\setminus F\imi$.  For $u\in 
F^0_i\setminus F^0\imi$,
$$\eqalign{\pi_n\circ\psi_n(u)&=\pi_n\bigl(c_\hi(u)\bigr)\cr
&=\pi_n(w_uw_u^*)\cr
&=r(u)r(u)^*\cr
&=\theta_n\bigl(\phi\nplpl(u)\bigr).\cr}$$
For $e\in F^1_i\setminus F^1\imi$ with $o(e)=u\in D$,
$$\eqalign{\pi_n\circ\psi_n(ee^*)&= 
\pi_n\bigl(c_\hi(e)c_\hi(e)^*\bigr)\cr
&=\pi_n\bigl(\eta(e)\eta(e)^*\bigr)\cr
&=\pi_n\bigl(w_e\psi\nmi(e)w_{t(e)}^*\bigr)
\pi_n\bigl(w_e\psi\nmi(e)w_{t(e)}^*\bigr)^*\cr
&=s(e)\pi_n\bigl(\psi\nmi(e)\bigr) r\ofte^*r\ofte\psi\nmi(e)^* 
s(e)^*\cr
&=s(e)\pi_n\bigl(\psi\nmi(e)\bigr)\pi_n\bigl(\psi\nmi(e)\bigr)^*s(e)^*\cr
&=\theta_n\bigl(\phi\nplpl(ee^*)\bigr).\cr}$$

Finally, we note that for $u\in F^0_i\setminus F^0\imi$, the element 
$w_u$ satisfies
$$\eqalign{w_u^*w_u&\approx \psi\nmi(u)\cr
w_uw_u^*&=c_\hi(u)\cr
&\approx\psi_n(u).\cr}$$
Thus, by choosing a larger value for $\hi\in\Lambda$ in the last
application of Lemma 3.3,
if necessary, we may assume that 
$$\psi\nmi(u)\sim\psi_n(u),\quad u\in F^0_i\setminus F^0\imi.$$
By induction we obtain $\bigl\{\psi_n(e)\bigm|e\in E^1\bigr\}$ 
satisfying $(1_n)$ -- $(6_n)$ and $(\oh)$. \qed
\enddemo

\proclaim{Theorem 3.12}  If $A$ is a simple separable nuclear purely 
infinite \cstar-algebra satisfying the Universal Coefficient
Theorem, and $A$ has finitely generated K-theory and 
torsion-free $K_1$, then $A$ is semiprojective.
\endproclaim

\demo{Proof} We rely, of course, on the celebrated theorem of Kirchberg/Phillips
([K], [P]), stating that  simple, separable, nuclear, purely infinite 
\cstar-algebras with UCT that
have isomorphic K-theory are strongly Morita 
equivalent.  In addition, we have [BGR], which  
tells us that strongly Morita equivalent separable \cstar-algebras are 
stably isomorphic.

The case where $\rank K_0(A)=\rank K_1(A)$ is proved in [B1].  In 
the remaining cases, it follows from Theorem 3.1 and Theorem 2.3 that 
there is a nonunital \cstar-algebra strongly Morita equivalent to $A$ 
that is semiprojective.  Since nonunital simple purely infinite 
\cstar-algebras are stable (by a theorem of Zhang ([Z])), it follows 
that $\k\otimes A$ is semiprojective.  Thus we may assume that $A$ is 
unital.

Let $u_1$, $u_2$, $\ldots$ be isometries in $A$ with pairwise 
orthogonal final projections.  Then $A_0=\overline{\hbox{span}}\, 
\{u_iAu_j\mid1\le i,j\}$ is a hereditary subalgebra of $A$, and 
$A_0\cong\k\otimes A$.  Thus we know that $A_0$ is semiprojective.

Now let $B$ be a \cstar-algebra,
 let $\l$ be a directed 
family of ideals of $B$ with closure $I$, and suppose that $A=B/I$.  
Let $\pi:B\to B/I$ be the quotient map, and set $B_0=\pi\inv(A_0)$.  
Since $A_0$ is semiprojective there is $J\in\l$ and a 
$*$-homomorphism $\psi_0:A_0\to B_0/J$ with $\pi\circ\psi_0 = 
\hbox{id}_{A_0}$.  We use $\pi_J$ to denote the quotient map of $B$ 
onto $B/J$.  Let $v\in B$ with $\pi(v)=u_1$.
Increasing $J$ if necessary, we may assume that $\pi_J(v^*v)=1_{B/J}$ 
and $\pi(vv^*)=\psi_0(u_1u_1^*)$.  Now define $\psi:A\to B/J$ by
$$\psi(x)=v^*\psi_0(u_1xu_1^*)v.$$
Then $\psi$ is a $*$-homomorphism, and for $x\in A$,
$$\eqalign{\pi\circ\psi(x)&=\pi(v^*)\pi\circ\psi_0(u_1xu_1^*)\pi(v)\cr
&=u_1^*(u_1xu_1^*)u_1\cr
&=x.\cr}$$
Thus $A$ is semiprojective. (See also [N], Proposition 5.48.)
\qed
\enddemo

\head References \endhead

\item{[B1]} B. Blackadar, Shape theory for \cstar-algebras, {\it Math. 
Scand.\/} {\bf 56} (1985), 249-275.
\item{[B2]} B. B.ackadar, Semiprojectivity in simple  
\cstar-algebras, to appear in the Proceedings of the U.S.-Japan Seminar.
\item{[BGR]} L.G. Brown, P. Green and M. Rieffel, Stable isomorphism 
and strong Morita equivalence of \cstar-algebras, {\it Pacific J. Math.\/} 
{\bf 71} (1977), 349-363.
\item{[C1]} J. Cuntz, A class of \cstar-algebras and topological Markov 
chains II:  Reducible chains and the Ext-functor for \cstar-algebras, 
{\it Invent. Math.\/} {\bf 63} (1981), 25-40.
\item{[C2]} J. Cuntz, K-theory for certain \cstar-algebras, {\it Ann. 
Math.\/} {\bf 113} (1981), 181-197.
\item{[CK]} J. Cuntz and W. Krieger, A class of \cstar-algebras and 
topological Markov chains, {\it Invent. Math.\/} {\bf 56} (1980), 
251-268.
\item{[EK]} E.G. Effros and J. Kaminker, Homotopy continuity and 
shape theory for \cstar-algebras, in {\it Geometric Methods in Operator 
Algebras\/}, eds. Araki and Effros, Pitman Res. Notes Math. {\bf 123},
Longman, Harlow, 1986.
\item{[EL]}, R. Exel and M. Laca, Cuntz-Krieger algebras for infinite 
matrices, {\it J. reine u. angew. Mathematik\/} {\bf 512}, (1999), 
119-172.
\item{[G]} J. Glimm, On a certain class of operator algebras, {\it 
Trans. Amer. Math. Soc.\/} {\bf 95} (1960), 318-340.
\item{[K]} E. Kirchberg, The classification of purely infinite 
\cstar-algebras using Kasparov's theory, (preprint).
\item{[N]} B. Neub\"{u}ser, Semiprojektivit\"{a}t und realisierunen 
von rein unendlichen \cstar-alge\-bren, preprint, M\"{u}nster, 2000.
\item{[P]} N.C. Phillips, A classification theorem for nuclear purely 
infinite simple \cstar-algebras, {\it Documenta Math.\/} {\bf 5} (2000), 
49-114.
\item{[RS]} I. Raeburn and W. Szymanski, Cuntz-Krieger algebras of 
infinite graphs and matrices, preprint.
\item{[R]} M. R\o rdam, Classification of Cuntz-Krieger algebras, {\it 
K-Theory\/} {\bf 9} (1995), 31-58.
\item{[Sp]} J. Spielberg, A functorial approach to the \cstar-algebras 
of a graph, {\it Internat. J. Math.\/} {\bf 13} 
(2002), no.3, 245-277.
\item{[Sz1]} W. Szymanski, On semiprojectivity of \cstar-algebras of 
directed graphs, Proc. AMS (to appear).
\item{[Sz2]} W. Szymanski, The range of $K$-invariants for 
\cstar-algebras of infinite graphs, preprint.
\item{[Z]} S. Zhang, Certain \cstar-algebras with real rank zero and 
their corona and multiplier algebras, {\it Pacific J. Math.\/} {\bf 
155} (1992), 169-197.
\enddocument
\bye